\documentclass[11pt]{preprint}
\usepackage{ScriptFonts}
\usepackage{times}
\usepackage{mhequ}
\usepackage{mhenvs}
\usepackage{mhfig}
\usepackage{Symbols}

\edef\myout{%
	\setbox255\vbox{%
	\unvbox255}%
	\the\output}

\def\vec#1{{\stackrel{\vbox{\hbox{\kern0.04em\tiny$\rightarrow$}\kern-0.1em}}{#1}}}
\def\cev#1{{\stackrel{\vbox{\hbox{\kern-0.04em\tiny$\leftarrow$}\kern-0.1em}}{#1}}}

\newcounter{assum}
\let\oldassum\theassum
\def\theassum{{\bf A\oldassum}}
\newenvironment{assumption}
{\refstepcounter{assum}%
\list{}{
	\setlength{\leftmargin}{\parindent}
	\setlength{\labelwidth}{\parindent}
	\setlength{\labelsep}{0mm}
	}
\item[{\theassum}\hfill]\it}
{\endlist}

\newcounter{bassum}
\let\oldbassum\thebassum
\def\thebassum{{\bf B\oldbassum}}
\newenvironment{bassumption}
{\refstepcounter{bassum}%
\list{}{
	\setlength{\leftmargin}{\parindent}
	\setlength{\labelwidth}{\parindent}
	\setlength{\labelsep}{0mm}
	}
\item[{\thebassum}\hfill]\it}
{\endlist}

\def\dstar{{\vbox{\kern0.12em\scriptsize\hbox{\kern-0.1em$\star$}}}}
\def\HS{{\hbox{\scriptsize{\rm HS}}}}
\def\XX{\mathsf X}
\def\YY{\mathsf Y}
\def\sC{Q}
\let\PH\Phi
\def\SP{{\Omega}}
\def\FF{{\cscr F}}
\def\PP{{\mathsf P}}
\def\WW{{\mathsf W}}
\def\pP^#1_#2{{\mathsf P}^{#1}_{\!#2}}
\def\pPt_#1{{\mathsf P}_{\!#1}^{}}
\def\hPt{\hat{{\mathbf P}}}
\def\pQ{{\bf Q}}
\def\pR{{\bf R}}
\def\nP^#1_#2{{\mathbf P}^{#1}_{\!#2}}
\def\Exp{{\bf E}}
\def\Meas{{\cscr M}}
\def\from{\,\colon}
\def\PC#1{{\bf C}^\infty_{#1}}
\def\pC{{\bf C}}
\let\Dens\CD
\def\DD#1#2{\Dens_{#1,#2}}
\def\Lip{{\rm L}^{}}
\def\Dom{{\cscr D}}
\def\GN{{G_{\!N}}}

\def\dmap#1#2{\psi_{#1\hbox{\tiny$\to$}#2}}
\def\imap#1#2{\psi_{#1\hbox{\tiny$\leftarrow$}#2}}
\def\dPsi#1#2{\Psi_{#1\hbox{\tiny$\to$}#2}}
\def\iPsi#1#2{\Psi_{#1\hbox{\tiny$\leftarrow$}#2}}
\def\dW#1#2{W_{#1\hbox{\tiny$\to$}#2}}
\def\iW#1#2{W_{#1\hbox{\tiny$\leftarrow$}#2}}

\begin{document}
\title{Exponential Mixing Properties of Stochastic PDEs\\Through Asymptotic Coupling}
\author{M.~Hairer}
\date{September 14, 2001}
\institute{\footnotesize D\'epartement de Physique Th\'eorique\\
Universit\'e de Gen\`eve, 1211 Gen\`eve 4, Switzerland\\
E-mail: Martin.Hairer@physics.unige.ch}

\maketitle

\begin{abstract}
We consider parabolic stochastic partial differential equations driven
by white noise in time. We prove exponential convergence of the transition probabilities towards a unique invariant measure under suitable conditions. These conditions amount essentially to the fact that the equation transmits the noise to all its determining modes. Several examples are investigated, including some where the noise does {\it not} act on every determining mode directly.
\end{abstract}

\tableofcontents

\thispagestyle{empty}
\newpage

\section{Introduction}
\label{sect:intro}

We are interested in the study of long-time asymptotics for parabolic stochastic partial differential equations. More precisely, the existence, uniqueness, and speed of convergence towards the invariant measure for such systems is investigated. The general setting is that of a stochastic PDE of the form
\begin{equ}[e:basic]
 dx = Ax\,dt + F(x)\,dt + Q\,d\omega(t)\;,\qquad x(0) = x_0\;,
\end{equ}
where $x$ belongs to some Hilbert space $\CH$, $A$ is the generator of a $C_0$-semigroup on $\CH$, $F\from\CH\to\CH$ is some nonlinearity, $\omega $ is the cylindrical Wiener process on some other Hilbert space $\CW$, and $Q\from\CW\to\CH$ is a bounded operator. If the nonlinearity $F$ is sufficiently ``nice'', there exists a unique solution $x(t)$ to \eref{e:basic} (see \eg \cite{ZDP1}).
In this paper, we investigate the asymptotic stability of \eref{e:basic}. We say that the solutions of \eref{e:basic} are asymptotically stable if there exists a {\it unique} probability measure $\mu_*$ on $\CH$ such that the laws of $x(t)$ converge to $\mu_*$, independently of the initial condition $x_0$.
We are interested in the situation where the asymptotic stability is a consequence of the noise (\ie the deterministic equation $\dot x = A x + F(x)$ is not asymptotically stable in the above sense), although the
noise is weak, in the sense that the range of $Q$ in $\CH$ is ``small''.

The investigation of asymptotic behaviour for solutions of \eref{e:basic} goes back to the early eighties (see for example \cite{Masl} for an excellent review article or the monograph \cite{ZDP} for a detailed exposition). Until recently, two approaches dominated the literature on this subject. For the first approach, sometimes called the ``dissipativity method'', one considers two solutions $x(t)$ and $y(t)$ of \eref{e:basic}, corresponding to the same realization of the Wiener process $\omega$, but with different initial conditions $x_0$ and $y_0$. If $A$ and $F$ are sufficiently dissipative, $\|x(t)-y(t)\|$ converges to $0$ for large times in some suitable sense. If this convergence is sufficiently fast and uniform, it yields asymptotic stability results (see for example \cite{DZA}). Closely related to this approach are the Lyapunov function techniques, developed for semilinear equations in \cite{Ich}. The dissipativity method, as well as the Lyapunov function techniques, are limited by the requirement that the deterministic equation $\dot x = Ax + F(x)$ already shows stable behaviour.

The (linearly) unstable situations are covered by the second approach, to which we refer as the ``overlap method''. It consists in showing that the Markov transition semigroup associated to \eref{e:basic} has the strong Feller property and is topologically irreducible. Then, provided that the equation \eref{e:basic} shows some dissipativity, arguments as developed in the monograph \cite{MT}, allow to bound the overlap between transition probabilities starting at two different initial points. This in turn yields strong asymptotic stability properties.
The main technical difficulty of this approach is to show that the strong Feller property holds. This difficulty is usually mastered either by studying the infinite-dimensional backward Kolmogorov equation associated to \eref{e:basic} \cite{DZ91}, or by showing that the Markov transition semigroup has good smoothing properties \cite{DEZ,Cerr}. This technique is limited by the requirement that the noise be sufficiently non-degenerate. A typical requirement is that the range of $Q$ {\it contains} the domain of some positive power of $-A$. To our knowledge, only one work \cite{EH3,H01} shows the strong Feller property for a stochastic PDE in a situation where the range of $Q$ is not dense in $\CH$ (but still of finite codimension).

Very recently, a third approach, to which we refer as the ``coupling method'', emerged in a series of papers on the 2D Navier-Stokes equation. (See \cite{ArmenKuk,MatNS,MY} and the references in Section~\ref{sect:examples}.) The main idea of these papers is to make a splitting $\CH = \CH_L\oplus\CH_H$ of the dynamics into a finite-dimensional, linearly unstable, low-frequency part $\CH_L$ and a remaining infinite-dimensional stable part $\CH_H$. An important assumption on $Q$ is then that the range of $Q$ contains $\CH_L$.
The space $\CH_L$ is chosen in such a way that the long-time asymptotics of the dynamics is completely dominated by the behaviour of the low-frequency part. More precisely, for any given realization $x_L(t)$ of the low-frequency part, the dynamics of the high-frequency part $x_H(t)$ will loose memory of its initial condition exponentially fast. On the low-frequency part, in turn, the noise acts in a non-degenerate way. A clever coupling argument allows to combine these two facts in order to obtain asymptotic stability results. The argument consists in coupling two realizations of \eref{e:basic} in such a way that if the low-frequency parts meet at some time $\tau$, they remain equal for all times $t>\tau$. (Of course, one has to show that $\tau$ is finite with probability $1$.) In fact, this coupling method is very close to the Gibbsian approach developed in \cite{KS00,BKLExp,EMS}, which consisted in transforming the infinite-dimensional Markovian system on $\CH$ to a finite-dimensional non-Markovian system on $\CH_L$. This finite-dimensional system was shown to have exponentially decaying memory and thus techniques from statistical mechanics can be applied.

Loosely speaking, the coupling method combines the arguments of both the dissipativity method (on $\CH_H$) and the overlap method (on $\CH_L$).
The coupling method thus yields a very powerful approach to the problem of asymptotic stability of \eref{e:basic}. The conditions of applicability of this coupling method have been successively weakened in the aforementioned papers, but the existing results always require, as we already mentioned, that the noise acts directly and independently on {\it every} determining mode of the equation. In this paper, we extend the coupling method to problems which do not satisfy this condition. Our overall approach is similar to the one exposed by Mattingly in \cite{MatNS}, and consequently some of our proofs are closely related to the arguments exposed there. Our main new idea is to construct a coupling for which the low-frequency parts of the dynamics do not actually meet at some finite time, but converge exponentially fast towards each other. This ``asymptotic coupling'' is achieved through a binding construction exposed in \sect{sect:couplconstr}, which seems to be new and can in some cases be implemented even in very degenerate situations.

In the following section, we illustrate the method of asymptotic coupling for a simple finite dimensional problem.

\subsection{A toy model}
\label{sect:toy}

Consider the following
system of stochastic differential equations in $\R^2$:
\begin{equa}[e:toy]
 dx_1 &= (2x_1 + x_2 - x_1^3)\,dt + d\omega(t)\;,\\
 dx_2 &= (2x_2 + x_1 - x_2^3)\,dt\;.
\end{equa}
This equation should be interpreted in the integral sense, with $\omega\in\Omega$ a Brownian motion.
Applying H\"ormander's condition \cite{Ho,Norr}, it is easy to see that the transition probabilities of \eref{e:toy} are smooth with respect to the Lebesgue measure on $\R^2$. Furthermore, an easy controllability argument shows that they have support everywhere and therefore are all mutually equivalent. Since \eref{e:toy} also exhibits a strong drift towards the center of
the phase space at large amplitudes, it follows by standard arguments that \eref{e:toy} possesses a unique invariant measure $\mu_*$ and that every initial condition is exponentially (in variation norm) attracted by $\mu_*$.

The problem with this argument is that it heavily relies on the existence of some reference measure (in this case the Lebesgue measure) which is equivalent to the transition probabilities. In the infinite-dimensional setting, such a reference measure will usually not exist when the noise is sufficiently degenerate. (For an account of some cases where such a reference measure does exist in the infinite-dimensional case, see \cite{Masl,EH3}.)
Furthermore, the fact that both directions in \eref{e:toy} are linearly 
unstable prevents one from applying the coupling method as it is
presented in the previous section.

We will show that the invariant measure for \eref{e:toy} is unique, using a coupling construction which pushes solutions together at an exponential rate. This construction is asymptotic, compared to more conventional coupling constructions, which look for hitting times at which the coupled dynamics actually meets.

Before we proceed, let us explain briefly what we mean by ``coupling''. A coupling for \eref{e:toy} is a process $(x(t),y(t))\in\R^2\times\R^2$, whose marginals $x(t)$ and $y(t)$ taken separately are both solutions of \eref{e:toy} (but with different initial conditions). In general, one takes a measure $\PP$ on $\Omega\times\Omega$, whose marginals are both equal to the Wiener measure $\WW$. Then a coupling for \eref{e:toy} can be constructed by drawing a pair $(\omega,\tilde\omega)\in\Omega\times\Omega$ distributed according to $\PP$ and solving the equations
\begin{equs}[2]
 dx_1 &= (2x_1 + x_2 - x_1^3)\,dt + d\omega(t)\;,&\qquad
 dy_1 &= (2y_1 + y_2 - y_1^3)\,dt + d\tilde\omega(t)\;,\\
 dx_2 &= (2x_2 + x_1 - x_2^3)\,dt\;,&
 dy_2 &= (2y_2 + y_1 - y_2^3)\,dt\;.\label{e:toycoupl1}
\end{equs}
We will carefully choose the measure $\PP$ in such a way that the quantity $\|x-y\|$ converges exponentially to $0$ for large times. This then yields the uniqueness of the invariant measure for \eref{e:toy}.

Our main idea leading to the construction of $\PP$ is to consider the following system in $\R^4$:
\begin{equa}[e:toycoupl]
 dx_1 &= (2x_1 + x_2 - x_1^3)\,dt + d\omega(t)\;,\\
 dx_2 &= (2x_2 + x_1 - x_2^3)\,dt\;,\\
 dy_1 &= (2y_1 + y_2 - y_1^3)\,dt + d\omega(t) + G(x_1,x_2,y_1,y_2)\,dt\;,\\
 dy_2 &= (2y_2 + y_1 - y_2^3)\,dt\;,
\end{equa}
where $d\omega$ denotes twice the same realization of the Wiener process. 
We see that this equation is the same as \eref{e:toycoupl1} with $\tilde \omega$ defined by
\begin{equ}[e:defwtilde]
 \tilde \omega(t) = \omega(t) + \int_0^t G\bigl(x_1(s),x_2(s),y_1(s),y_2(s)\bigr)\,ds\;.
\end{equ}
The noise $\tilde\omega\in\Omega$ is distributed according to some measure $\tilde\WW$ which is in general not equal to the Wiener measure $\WW$. Therefore, \eref{e:toycoupl} does not yet define a coupling for \eref{e:toy}. If $G$ is small in the sense that the quantity
\begin{equ}[e:GL2]
 \int_0^\infty \bigl\|G\bigl(x_1(s),x_2(s),y_1(s),y_2(s)\bigr)\bigr\|^2\,ds
\end{equ}
is bounded with sufficiently high probability, then the measures $\tilde\WW$ and $\WW$ are equivalent. In this case, it is possible to construct a measure $\PP$ on $\Omega\times\Omega$ whose marginals are $\WW$, with the important property that there exists a random time $\tau$ with ${\bf P}(\tau < \infty) = 1$ such that the solutions of the coupled system satisfy \eref{e:toycoupl} for times $t \ge \tau$.

In view of the above, we have reduced the problem to finding a function $G$ such that the solutions of \eref{e:toycoupl} satisfy $\|y(t)-x(t)\| \to 0$ for $t\to\infty$ and \eref{e:GL2} is bounded. We introduce the difference process $\rho = y - x$, and we write
\minilab{e:diff}
\begin{equs}
 \dot\rho_1 &= 2\rho_1 + \rho_2 - \rho_1\bigl(x_1^2 + x_1y_1 + y_1^2\bigr) + G(x,y)\;,\label{e:diff1}\\
 \dot\rho_2 &= 2\rho_2 + \rho_1 - \rho_2\bigl(x_2^2 + x_2y_2 + y_2^2\bigr)\;.\label{e:diff2}
\end{equs}
It is easy to find a function $G$ such that $\rho_1 \to 0$, but this does not yet mean that $\rho_2$ will go to zero. A closer look at \eref{e:diff2} shows that if we could force $\rho_1$ to be very close to $-3\rho_2$, \eref{e:diff2} could be written as
\begin{equ}
  \dot\rho_2 = -\rho_2 + \eps - \rho_2\bigl(x_2^2 + x_2y_2 + y_2^2\bigr)\;,
\end{equ}
which is asymptotically stable. Introduce the function $\zeta = \rho_1 + 3\rho_2$. We then have
\begin{equ}
 \dot \zeta = (\ldots) + G(x_1,x_2,y_1,y_2)\;,
\end{equ}
with $(\ldots)$ an expression of the order $\|\rho\|\bigl(1+\|x\|^2 + \|y\|^2\bigr)$. Now we can of course choose $G = - (\ldots) - 2\zeta$. This way, the equation for $\zeta$ becomes $\dot\zeta = -2\zeta$ and we have the solution $\zeta(t) = \zeta(0) e^{-2t}$. Plugging this into \eref{e:diff2}, we get
\begin{equ}
  \dot\rho_2 = -\rho_2 + \zeta(0)e^{-2t} - \rho_2\bigl(x_2^2 + x_2y_2 + y_2^2\bigr)\;.
\end{equ}
We thus have the estimate
\begin{equ}
 |\rho_2(t)| \le |\rho_2(0)|e^{-t} + |\zeta(0)| e^{-2t}\;.
\end{equ}
Finally, $\rho_1$ is estimated by using the definition of $\zeta$ and we get
\begin{equ}
 |\rho_1(t)| \le |\rho_2(0)|e^{-t} + 4|\zeta(0)| e^{-2t}\;.
\end{equ}
This shows that, with $G$ chosen this way, there exists a constant $C$ such that
\begin{equ}
 \|x(t) - y(t)\| \le C\|x(0) - y(0)\|e^{-t}\;,
\end{equ}
for almost every realization of the noise. Since typical realizations of $x(t)$ do not grow faster than linearly, $G$ is also of the order $e^{-t}$, with at most a polynomial factor in $t$ multiplying the exponential. The main result of this paper, \theo{theo:main}, shows that the above construction implies the existence and uniqueness of an invariant probability measure $\mu_*$ for the problem at hand. Furthermore, it shows that the transition probabilities converge exponentially fast towards $\mu_*$ in the Wasserstein norm (the dual norm to the Lipschitz norm on functions).

This concludes our presentation of the toy model. For a more precise statement, the reader is encouraged to actually check that the above construction allows to verify the assumptions stated in \sect{sect:equadiff}.

The remainder of this paper is organized as follows. In \sect{sect:coupl}, we give the precise definitions for the type of coupling we will consider. In \sect{sect:Prop}, we state the properties of the coupling that are required for our purpose. In \sect{sect:abstr}, we prove the abstract formulation of our main ergodic theorem. In \sect{sect:equadiff}, this abstract theorem is then specialized to the case of stochastic differential equations. In \sect{sect:examples} finally, we present several examples where our construction applies, although the noise does not act directly on every determining mode of the equation.

\begin{acknowledge}
The author thanks Jean-Pierre Eckmann, Guillaume van Baalen, Emmanuel Zabey, and Jacques Rougemont for several helpful discussions and comments. This work was partially supported by the Fonds National Suisse.
\end{acknowledge}

\section{The Coupling Construction}
\label{sect:coupl}

In this section, we explain our coupling construction. Before we start with the actual definitions of the various objects appearing in the construction, we fix our notations.

\subsection{Notations}\label{sect:not}

If $\mu$ is a measure on a measurable space $\XX$ (in the sequel, we will always consider Polish\footnote{\ie complete, separable, and metric} spaces) and $f\from\XX \to \YY$ is a measurable map, we denote by $f^*\mu$ the measure on $\YY$ defined by $(f^*\mu)(A) \equiv \mu(f^{-1}(A))$. For example, if $\Pi$ is a projection on one component of a product space, $\Pi^*\mu$ denotes the marginal of $\mu$ on this component. If a natural reference measure is given on the current space, we denote by $\Dens\mu$ the density of $\mu$ with respect to the reference measure.

We define for any two measures $\mu$ and $\nu$ the measures $\mu\wedge\nu$ and $\mu \setminus \nu$. If a common reference measure is given, these operations act on densities like
\begin{equs}
 \bigl(\Dens(\mu\wedge\nu)\bigr)(x) &= \min\{\Dens\mu(x)\;,\;\Dens\nu(x)\}\;,\\
 \bigl(\Dens(\mu\setminus\nu)\bigr)(x) &= \max\{\Dens\mu(x)-\Dens\nu(x),0\}\;.
\end{equs}
It immediately follows that $\mu = (\mu\wedge\nu) + (\mu\setminus\nu)$ for any two measures $\mu$ and $\nu$.
We will use the equivalent notations $\mu \le \nu$ and $\nu \ge \mu$ to say that $\mu \wedge \nu = \mu$ holds. One can check the following relations:
\begin{equa}
 f^*(\mu \wedge \nu) \le f^*\mu \wedge f^*\nu\;,\\
f^*(\mu \setminus\nu) \ge f^*\mu \setminus f^*\nu\;.
\end{equa}
Equalities hold if $f$ is injective.

For a given topological space $\XX$, we denote by $\Meas(\XX)$ the space of all finite signed Borel measures on $\XX$. We denote by $\Meas_1(\XX)$ the set of all probability measures on $\XX$. For $\mu\in\Meas(\XX)$, we denote by $\|\mu\|$ its total variation norm (which is simply its mass if $\mu$ has a sign).

\subsection{Definition of coupling}

In this section, and until the end of the paper, we will often consider families $\pQ_y$ of measures indexed by elements $y\in\YY$, with $\YY$ some Polish space. One should think of $y$ as the initial condition of a Markov chain on $\YY$ and of $\pQ_y$ either as its transition probabilities, or as the measure on pathspace obtained by starting from $y$. We will always assume that the functions $y\mapsto \pQ_y(A)$ are measurable for every Borel set $A$. If $\pQ_y$ is a family of measures on $\YY^n$ and $\pR_y$ is a family of measures on $\YY^m$, a family of measures $(\pR\pQ)_y$ on $\YY^{n+m} = \YY^n\times\YY^m$ can be defined on cylindrical sets in a natural way by
\begin{equ}[e:comp]
 \bigl(\pR\pQ\bigr)_y(A\times B) = \int_{A}\pR_{z_n}(B)\,\pQ_y(dz)\;,
\end{equ}
where $A \subset \YY^n$, $B \subset \YY^m$, and $z_n$ denotes the $n$th component of $z$.

We consider a discrete-time Markovian random dynamical system (RDS) $\PH$ on a Polish space $\XX$ with the following structure. There exists a ``one-step'' probability space $(\SP,\FF,\PP)$ and $\PH$ is considered as a jointly measurable map $\PH \from (\XX,\SP) \to \XX$. The iterated maps $\PH^n \from (\XX, \SP^n) \to \XX$ with $n \in \N$ are constructed recursively by 
\begin{equ}
 \PH^n(x,\omega_1,\ldots,\omega_n) = \PH\bigl(\PH^{n-1}(x,\omega_1,\ldots,\omega_{n-1}),\omega_n\bigr)\;,
\end{equ}
This construction gives rise to a Markov chain on $\XX$ (also denoted by $\PH$) with one-step transition probabilities
\begin{equ}
 \pP^_x \equiv \PH(x,\cdot\,)^*\PP\;.
\end{equ}
The $n$-step transition probabilities will be denoted by $\pP^n_x$.
Our main object of study will be the family of measures on pathspace generated by $\PH$. Take a sequence $\{\omega_i\}_{i=0}^\infty$ and an initial condition $x\in\XX$. We then define $x_0 = x$ and $x_{i+1} = \PH(x_i,\omega_i)$. 
 We will denote by $\nP^n_x$ with $n\in\N\cup \{\infty\}$ the measure on $\XX^n$ obtained by transporting $\PP^n$ with the map $\{\omega_i\} \mapsto \{x_i\}$. 
It is also natural to view $\nP^n_x$ as a measure on $\XX^n \times \SP^n$ by transporting $\PP^n$ with the map $\{\omega_i\} \mapsto \{x_i,\omega_i\}$, so we will use both interpretations.

\begin{remark}
The above setup is designed for the study of stochastic differential equations driven by additive noise. In that case, $\SP$ is some Wiener space and $\PH$ maps an initial condition and a realization of the Wiener process on the solution after time $1$. Nevertheless, our setup covers much more general cases.
\end{remark}

The coupling needs two copies of the pathspace, \ie we will consider elements $(x,y) \in \XX^\infty \times \XX^\infty$. It will be convenient to use several projectors from $\XX^N \times \XX^N$ to its components. We define therefore (for $n \le N$):
\begin{equ}
 \Pi_1 \from (x,y)\mapsto x\;,\quad \Pi_2 \from (x,y)\mapsto y\;,\quad
\pi_n \from (x,y) \mapsto (x_n,y_n)\;.
\end{equ}
We also define $\pi_{i,n} \equiv \Pi_i\circ\pi_n$ for $i\in \{1,2\}$.

\begin{definition}
Let $\PH$ be a Markov chain on a Polish space $\XX$ and let $\nP^{\infty}_{x}$ be the associated family of measures on the pathspace $\XX^\infty$. A {\em coupling} for $\PH$ is a family $\PC{x,y}$ of probability measures on $\XX^\infty \times \XX^\infty$ satisfying
\begin{equ}
 \Pi_1^*\PC{x,y} = \nP^{\infty}_{x}\qquad\text{and}\qquad
 \Pi_2^*\PC{x,y} = \nP^{\infty}_{y}\;,
\end{equ}
where $\Pi_1$ and $\Pi_2$ are defined as above.
\end{definition}

A trivial example of coupling is given by $\PC{x,y} = \nP^{\infty}_{x}\times \nP^{\infty}_{y}$.
The interest of constructing a non-trivial coupling comes from the following observation.
Take some suitable set of test functions $\CG$ on $\XX$ and define a norm on $\Meas(\XX)$ by
\begin{equ}
 \|\mu\|_{\CG} = \sup_{g \in \CG} \,\scal{g,\mu}\;.
\end{equ}
Once the existence of an invariant measure for the Markov chain $\PH$ is established, one usually wishes to show its uniqueness by proving that $\PH$ forgets about its past sufficiently fast, \ie
\begin{equ}
 \lim_{n\to \infty} \|\pP^n_x - \pP^n_y\|_{\CG} = 0\;,\quad\text{for all}\quad (x,y) \in \XX^2\;,
\end{equ}
with suitable bounds on the convergence rate as a function of the initial conditions.
Now take a coupling $\PC{x,y}$ for $\PH$. It is straightforward to see that by definition the equality
\begin{equ}
 \scal{\pP^n_x,g} = \int_{\XX\times\XX} g(z)\,\bigl(\pi_{1,n}^*\PC{x,y}\bigr)(dz)
\end{equ}
holds, as well as the same equality where $\pi_{1,n}$ is replaced by $\pi_{2,n}$ and $\pP^n_x$ is replaced by $\pP^n_y$. Therefore, one can write
\begin{equ}[e:couplingDiff]
\|\pP^n_x - \pP^n_y\|_{\CG}  = \sup_{g\in\CG}\int_{\XX\times\XX} \bigl(g(\Pi_1z)-g(\Pi_2z)\bigr)\,\bigl(\pi_{n}^*\PC{x,y}\bigr)(dz)\;.
\end{equ}
This equation is interesting, because it is in many cases possible to construct a coupling $\PC{x,y}$ such that for $n$ large, the measure $\pi_{n}^*\PC{x,y}$ is concentrated near the diagonal $\Pi_1 z = \Pi_2 z$, thus providing through \eref{e:couplingDiff} an estimate for the term $\|\pP^n_x - \pP^n_y\|_{\CG}$. This is precisely what was shown in our toy model of \sect{sect:toy}, where we constructed $f$ in such a way that $\|x(t)-y(t)\| \to 0$ for $t\to \infty$.

\subsection{The binding construction}
\label{sect:couplconstr}

In this subsection, we describe a specific type of coupling for a given RDS $\PH$. Only couplings of that type will be taken under consideration in the sequel.

Let $\PH$ and the associated probability space $(\SP,\FF,\PP)$ be as above. We consider a family $\dmap{x}{y}\from\SP\to\SP$ (the pair $(x,y)$ belongs to $\XX^2$) of measurable functions that also have measurable inverses. We will call these functions {\em binding functions} for $\PH$. The reason for this terminology is that, given a realization $\{\omega_n\}_{n=0}^\infty$ of the noise and a pair of initial conditions $(x_0,y_0) \in \XX^2$, the binding functions allow us to construct two paths $\{x_n\}$ and $\{y_n\}$ by setting
\begin{equ}[e:defCouply]
\tilde \omega_n =\dmap{x_n}{y_n}(\omega_n)\;,\qquad  x_{n+1} = \PH(x_{n},\omega_n)\;,\qquad y_{n+1} = \PH(y_{n},\tilde\omega_n)\;.
\end{equ}
Our aim is to find a family $\dmap{x}{y}$ such that $y_n$ converges towards $x_n$ in a suitable sense for large values of $n$. Thus, the binding functions play the role of a spring between $x$ and $y$.
We will say that \eref{e:defCouply} is a {\em binding construction} for $\PH$. We denote the inverse of $\dmap{x}{y}$ by $\imap{x}{y}$. The reason behind this notation should be clear from the diagram below.
\begin{equ}[e:diag]
 \mhpastefig{Diagramme}
\end{equ}
The solid arrows denote the various maps and the dashed arrows denote the influences of the appearing quantities on those maps.
It shows that it is also possible to achieve the binding construction by
first choosing a sequence $\{\tilde\omega_n\}_{n=0}^\infty$ and then using
$\imap{x_n}{y_n}$ to construct the $\omega_n$, thus obtaining the same set of possible realizations for $(x_n,y_n)$. This symmetry between $\dmap{x}{y}$ and $\imap{x}{y}$ is also explicit in \eref{e:defPxy} below.

Guided by the above construction, we use the binding maps to construct a coupling Markov chain $\Psi$ on $\XX\times\XX$ with transition probabilities $\pC_{x,y}$ in the following way. Define the maps
\begin{equa}[2][e:defPsi]
 \dPsi{x}{y} \from\SP &\to \SP\times\SP \qquad&\qquad\qquad
 \iPsi{x}{y} \from \SP &\to \SP\times\SP\\
\omega &\mapsto \bigl(\omega,\dmap{x}{y}(\omega)\bigr)\;,&
\omega&\mapsto \bigl(\imap{x}{y}(\omega),\omega\bigr)\;.
\end{equa}
Notice that, up to some null set, the image of both maps is the set $\{(\omega,\tilde\omega)\;|\; \tilde\omega = \dmap{x}{y}(\omega)\}$.
Then we define a family of measures $\pPt_{x,y}$ on $\SP\times\SP$ by
\begin{equ}[e:defPxy]
 \pPt_{x,y} = \bigl(\dPsi{x}{y}^*\PP\bigr) \wedge \bigl(\iPsi{x}{y}^*\PP\bigr) = \iPsi{x}{y}^*\bigl(\PP \wedge \dmap{x}{y}^*\PP\bigr)\;.
\end{equ}
According to \eref{e:diag}, the measure $\pPt_{x_n,y_n}$ is precisely the common part  between the measure obtained for $(\omega_n,\tilde\omega_n)$ by distributing $\omega_n$ according to $\PP$ and the one obtained by distributing $\tilde\omega_n$ according to $\PP$.
Thus both marginals of the measure $\pPt_{x,y}$ are smaller (in the sense of Section~\ref{sect:not}) than $\PP$.
In order to have a non-trivial construction, we impose that the measures $\PP$ and $\dmap{x}{y}^*\PP$ are equivalent. The density of $\dmap{x}{y}^*\PP$ relative to $\PP$ will be denoted by $\DD{x}{y}(\omega)$.

Considering again \eref{e:diag}, the family of measures $\pPt_{x,y}$ is transported on $\XX\times\XX$ by
defining 
\begin{equs}
 \Phi_{x,y}\from\SP\times\SP&\to\XX\times\XX \\
(\omega,\tilde\omega) &\mapsto \bigl(\PH(x,\omega),\PH(y,\tilde\omega)\bigr)\;,
\end{equs}
and setting
\begin{equ}[e:defCalt]
  \pQ_{x,y} \equiv \Phi_{x,y}^* \pPt_{x,y} \;.
\end{equ}
But this does not give a transition probability function yet, since the
measures $\pPt_{x,y}$ are not normalized to $1$. We therefore
define the family of measures $\nP^_{x,y}$ by
\begin{equ}
 \nP^_{x,y} = \pPt_{x,y} + c_{x,y}\, \bigl(\PP \setminus \Pi_1^*\pPt_{x,y}\bigr)
\times \bigl(\PP \setminus \Pi_2^*\pPt_{x,y}\bigr)\;,
\end{equ}
where the number $c_{x,y}$ is chosen in such a way that the resulting measure is a probability measure. By a slight abuse of notation, we used here the symbol $\Pi_i$ to denote the projection on the $i$th component of $\SP\times\SP$. As a matter of fact, $\bigl(\PP \setminus \Pi_1^*\pPt_{x,y}\bigr)$ and $\bigl(\PP \setminus \Pi_2^*\pPt_{x,y}\bigr)$ have the same mass, which is equal to $1-\|\pPt_{x,y}\|$, so
\begin{equ}
 c_{x,y} = {1\over \|\PP \setminus  \Pi_2^*\pPt_{x,y}\|}\;,
\end{equ}
for example. (Recall that the symbol $\|\cdot\|$ stands for the total variation norm, which is simply equal to its mass for a positive measure.)
It is straightforward to show that the following holds:

\begin{lemma}\label{lem:marg}
The measures $\nP^_{x,y}$ satisfy $\Pi_i^*\nP^_{x,y} = \PP$ for $i=1,2$.
\end{lemma}
\begin{proof}
It is clear by \eref{e:defPxy} that $\Pi_i^*\pPt_{x,y} \le \PP$. Thus
\begin{equa}[e:verif]
 \Pi_1^*\nP^_{x,y} &= \Pi_1^*\pPt_{x,y} + c_{x,y} \|\PP \setminus  \Pi_2^*\pPt_{x,y}\| \bigl(\PP \setminus  \Pi_1^*\pPt_{x,y}\bigr) \\
&= \bigl(\PP \wedge \Pi_1^*\pPt_{x,y}\bigr) + \bigl(\PP \setminus  \Pi_1^*\pPt_{x,y}\bigr) = \PP\;,
\end{equa}
and similarly for $\Pi_2^*\nP^_{x,y}$.
\end{proof}

This finally allows us to define the transition probabilities for $\Psi$ by
\begin{equ}[e:defCxy]
 \pC_{x,y} =  \Phi_{x,y}^* \nP^_{x,y} \equiv \pQ_{x,y} + \pR_{x,y}\;.
\end{equ}
In this expression, the only feature of $\pR_{x,y}$ we will use is that it is a positive measure.
We define $\pC^\infty_{x,y}$ as the measure on the pathspace $\XX^\infty\times\XX^\infty$ obtained by iterating \eref{e:comp}. Since $\Pi_1\circ\PH_{x,y} = \PH(x,\cdot)\circ\Pi_1$ and similarly for $\Pi_2$, it is straightforward to verify, using \lem{lem:marg}, that the measure $\pC^\infty_{x,y}$ constructed this way is indeed a coupling for $\PH$.

For a given step of $\Psi$, we say that the trajectories do couple if the step is drawn according to $\pQ_{x,y}$ and that they don't couple otherwise.

\begin{remark}\label{rem:prob}
Since $\nP^_{x,y}$ is a family of measures on $\SP\times\SP$, it is also possible to interpret $\pC^n_{x,y}$ as a family of probability measures on $\XX^n\times\XX^n \times \SP^n\times\SP^n$. We will sometimes use this viewpoint in the following section. It is especially useful when the RDS $\PH$ is obtained by sampling a continuous-time process.
\end{remark}

\begin{remark}\label{rem:Psihat}
It will sometimes be useful to have an explicit way of telling whether a
step of $\Psi$ is taken according to $\pQ_{x,y}^\infty$ or according to 
$\pR_{x,y}^\infty$ (\ie whether the trajectories couple or not). To this end, we introduce a Markov chain $\hat\Psi$
on the augmented phase space $\XX\times\XX\times\{0,1\}$ with transition probabilities
\begin{equ}
 \nP^_{x,y} = \pQ_{x,y}\times\delta_1 + \pR_{x,y}\times\delta_0\;.
\end{equ}
The marginal of $\hat\Psi$ on $\XX\times\XX$ is of course equal to $\Psi$.
By a slight abuse of notation, we will also write $\pC_{x,y}^\infty$ for
the probability measure on pathspace induced by $\hat\Psi$.
\end{remark}

It will be useful in the sequel to have a map that ``transports'' the family
of maps $\dmap{x}{y}$ on $\SP^n$ via the RDS $\PH$. More precisely, fix a pair $(x,y) \in \XX\times\XX$ of starting points and a sequence $(\omega_0,\ldots,\omega_n)$ of realizations of the noise. We then define $x_0=x$, $y_0 = y$, and, recursively for $i=0,\ldots,n$
\begin{equ}
 x_{i+1} = \PH(x_i,\omega_i)\;,\qquad y_{i+1} = \PH\bigl(y_i,\dmap{x_i}{y_i}(\omega_i)\bigr)\;.
\end{equ}
This allows us to define the family of maps $\Xi^n_{x,y} \from \SP^n \to \SP^n$ by
\begin{equ}[e:defXi]
	\Xi^{n+1}_{x,y}(\omega_0,\ldots,\omega_n) \mapsto \bigl(\dmap{x_0}{y_0}(\omega_0),\ldots,\dmap{x_n}{y_n}(\omega_n)\bigr)\;.
\end{equ}
Since $\dmap{x}{y}^*\PP$ is equivalent to $\PP$, we see that $\bigl(\Xi^n_{x,y}\bigr)^*\PP^n$ is equivalent to $\PP^n$ and we denote its density by $\DD{x}{y}^n$. We also notice that the family of measures $\pQ^n_{x,y}$ is obtained by transporting $\bigl(\Xi^n_{x,y}\bigr)^*\PP^n \wedge \PP^n$ onto $\XX^n\times\XX^n$ with the maps $\Phi_{x_i,y_i}\circ\dPsi{x_i}{y_i}$. In particular, one has the equality
\begin{equ}[e:relDens]
\|\pQ^n_{x,y}\| = \bigl\|\bigl(\Xi^n_{x,y}\bigr)^*\PP^n \wedge \PP^n\bigr\| = \int_{\SP^n} \bigl(1\wedge \Dens_{x,y}^n(\omega)\bigr)\,\PP^n(d\omega)\;.
\end{equ}

\section{Assumptions on the Coupling}
\label{sect:Prop}

In this section, we investigate the properties of the coupling $\pC_{x,y}^\infty$ constructed in the previous section. We give a set of assumptions on the binding functions $\dmap{x}{y}$ that ensure the existence and uniqueness of the invariant measure for $\PH$.

In order to achieve this, we want the map $\dmap{x}{y}$ to modify the noise in such a way that trajectories drawn according to $\pQ_{x,y}$ tend to come closer together. This will be the content of Assumption~\ref{ass:Conv}. Furthermore, we want to know that this actually happens, so the noise should not be modified too much. This will be the content of assumptions~\ref{ass:Prob} and \ref{ass:Coupl}. All these nice properties usually hold only in a ``good'' region of the phase space. Assumptions~\ref{ass:Lyap} and \ref{ass:W} will ensure that visits to this good region happen sufficiently often.

\subsection{Lyapunov structure}

Since we are interested in obtaining exponential mixing, we need assumptions of exponential nature. Our first assumption concerns the global aspects of the dynamics. It postulates that $\PH$ is attracted exponentially fast towards a ``good'' region of its state space. We achieve this by assuming the existence of a Lyapunov function for $\PH$.

\begin{definition}
Let $\PH$ by a RDS with state space $\XX$ as before. A {\em Lyapunov function} for $\PH$ is a function $V\from\XX\to[0,\infty]$ for which there exist constants $a\in(0,1)$ and $b>0$, such that
\begin{equ}[e:Lyap]
\int_{\SP}V\bigl(\PH(x,\omega)\bigr)\,\PP(d\omega) \le a V(x) + b\;,
\end{equ}
for every $x\in\XX$ with $V(x)<\infty$.
\end{definition}

Our first assumption then reads

\begin{assumption}\label{ass:Lyap}
There exist a Lyapunov function $V$ for $\PH$. Furthermore, $V$ is such that
\begin{equ}
\PP\bigl\{\omega\;|\;V\bigl(\PH(x,\omega)\bigr)<\infty\bigr\} = 1\;,
\end{equ}
for every $x\in\XX$.
\end{assumption}

For convenience, we also introduce the function $\tilde V \from \XX\times\XX\to[0,\infty]$ defined by
\begin{equ}
 \tilde V(x,y) = V(x) + V(y)\;.
\end{equ}
Notice that $\tilde V$ is a Lyapunov function for $\Psi$ by construction.

In some cases, when the control over the densities $\Dens_{x,y}$ is uniform enough or when the phase space is already compact (or bounded), a Lyapunov function is not needed. In such a situation, one can simply choose $V \equiv 1$.

In our case of interest, the RDS $\PH$ is obtained by sampling a continuous-time process $\Phi_t$ at discrete times. In that setting, it is useful to have means to control excursions to large amplitudes that take place between two successive sampling times. To this end, we introduce a function
$W\from\XX\times\SP\to[0,\infty]$ given by 
\begin{equ}
 W(x,\omega) = \sup_{t\in[0,1]} V\bigl(\Phi_t(x,\omega)\bigr)
\end{equ}
in the continuous-time setting and by
\begin{equ}
  W(x,\omega) = V(x)
\end{equ}
in the discrete-time setting. In fact, any other choice of $W$ is all right, as long as it satisfies the properties that are summarized in Assumption~\ref{ass:W} below.

Before stating these properties, we define two other functions that act on pairs
of initial conditions that couple by
\begin{equa}[e:defWxy]
\dW{x}{y}(\omega) &=  W(x,\omega) + W\bigl(y,\dmap{x}{y}(\omega)\bigr)\;,\\
\iW{x}{y}(\omega) &=  W\bigl(x,\imap{x}{y}(\omega)\bigr) + W(y,\omega)\;.
\end{equa}
We will assume that $W$ and the binding functions are such that $W$, $\dW{x}{y}$ and $\iW{x}{y}$ do not behave much worse than $V$. More precisely, we will assume that:

\begin{assumption}\label{ass:W}
There exists a function $W\from\XX\times\SP\to[0,\infty]$ such that
\minilab{e:propW}
\begin{equs}
\essinf_{\omega\in\SP}\,W(x,\omega) &= V(x)\;, \label{e:propW1}\\
\int_\SP W(x,\omega)\,\PP(d\omega) &\le c\,V(x)\;, \label{e:propW2}
\end{equs}
for some constant $c>0$.
Furthermore, there exist constants $C>0$ and $\delta \ge 1$ such that the estimates
\begin{equa}[e:propW3]
 \dW{x}{y}(\omega) &\le C\bigl(1 + V(y) + W(x,\omega)\bigr)^\delta\;,\\
 \iW{x}{y}(\omega) &\le C\bigl(1 + V(x) + W(y,\omega)\bigr)^\delta\;,
\end{equa}
hold for the functions defined in \eref{e:defWxy}.
\end{assumption}

The Lyapunov structure given by assumptions~\ref{ass:Lyap} and \ref{ass:W} ensures that $W$ (and thus also $V$) does not
increase too fast along a typical trajectory. In order to
make this statement precise, we define for a given initial condition $x\in\XX$ the sets $A_{x,k} \subset \SP^\infty$ by
\begin{equ}[e:defAxk]
 A_{x,k} = \bigl\{\omega \in \SP^\infty \;|\; W\bigl(\PH^n(x,\omega),\omega_n\bigr) \le k V(x) + kn^2\quad
\forall n > 0\bigr\}\;,
\end{equ}
where $k$ is some positive constant. The sets $A_{x,k}$ contain almost every typical reali\-zation of the noise:

\begin{lemma}\label{lem:bound}
Let $\PH$ be a RDS satisfying assumptions~\ref{ass:Lyap} and \ref{ass:W}. Then, there exists a constant $C>0$ such that
\begin{equ}
 \PP^\infty \bigl(A_{x,k}\bigr) \ge 1 - {C \over k}\;,
\end{equ}
for every $x\in\XX$ and every $k>0$.
\end{lemma}

\begin{proof}
For $\omega\in\SP^\infty$, we define $x_n = \PH^n(x,\omega)$.
 Notice that by \eref{e:propW2} and the Lyapunov structure, one has the estimate
\begin{equ}[e:estW]
 \Exp \bigl(W(x_n,\omega_{n+1})\bigr) \le c a^n V(x) + {bc \over 1-a}\;,
\end{equ}
where $\Exp$ denotes expectations with respect to $\PP^\infty$. We also notice that $A_{x,k} = \bigcap_{n > 0} A^{(n)}_{x,k}$ with
\begin{equ}
  A^{(n)}_{x,k} = \bigl\{\omega \;|\; W(x_n,\omega_{n+1}) \le	k V(x) + k n^2\bigr\}\;.
\end{equ}
Combining this with \eref{e:estW}, we see that 
\begin{equ}
 \PP^\infty \bigl( A^{(n)}_{x,k}\bigr) \ge 1 - {c\over k}{a^n V(x) + {b(1-a)^{-1}}\over V(x) + n^2}\;.
\end{equ}
Therefore, the worst possible estimate for $ \PP^\infty( A_{x,k})$ is
\begin{equ}
 \PP^\infty\bigl( A_{x,k}\bigr)  \ge 1 - {c\over k} \sum_{n=1}^\infty {a^n V(x) + {b(1-a)^{-1}}\over V(x) + n^2}\;,
\end{equ}
which proves the claim.
\end{proof}

\subsection{Binding property}

The crucial property of the coupling is to bring trajectories closer together. In order to make this statement more precise, we introduce the Lipschitz norm $\|\cdot\|_\Lip$ defined on functions $g\from\XX\to\R$ by
\begin{equ}
\|g\|_\Lip = \sup_{x\in\XX}|g(x)| + \sup_{x,y\in\XX}{|g(x)-g(y)|\over d(x,y)}\;,
\end{equ}
where $d(\cdot,\cdot)$ denotes the distance in $\XX$.
The dual norm on $\Meas(\XX)$ is then given by
\begin{equ}
\|\mu\|_\Lip = \sup_{\|g\|_\Lip = 1}\int_{\XX}g(x)\,\mu(dx)\;.
\end{equ}
With this definition at hand, we make the following assumption on
the coupling part $\pQ_{x,y}^\infty$.

\begin{assumption}\label{ass:Conv}
There exist a positive constant $\gamma_1$ and a family of constants $K \mapsto C_K$ such that, for every $K>0$,
\begin{equ}[e:convergence]
\|\pi_{1,n}^*\pQ_{x,y}^\infty - \pi_{2,n}^*\pQ_{x,y}^\infty\|_\Lip\le 
C_K e^{-\gamma_1 n}\;,
\end{equ}
holds when $\tilde V(x,y) \le K$.
\end{assumption}

\begin{remark}\label{rem:Conv}
The sub-probability kernels $\pQ_{x,y}$ are smaller than the transition probabilities for the binding construction \eref{e:defCouply}. Thus, \eref{e:convergence} is implied by an inequality of the type
\begin{equ}
 \Exp \bigl(d(x_n,y_n)\bigr) \le C\tilde V(x_0,y_0)e^{-\gamma_1 n}\;,
\end{equ}
where $d$ denotes the distance in $\XX$ and $\Exp$ denotes the expectation with respect to the construction \eref{e:defCouply}.
\end{remark}

Notice that this assumption is non-trivial only if our coupling is such that $\|\pQ_{x,y}^\infty\| > 0$ for sufficiently many starting points. This will be ensured by the next assumption.

\begin{assumption}\label{ass:Prob}
Let $\DD{x}{y}^n$ be defined as in Section~\ref{sect:couplconstr}. We assume that for every $K>0$,
there exists a family of sets $\Gamma^K_{x,y}\subset \SP^\infty$ and constants $c_1,c_2>0$ such that the estimates
\begin{equ}[e:assProb]
\PP^\infty(\Gamma^K_{x,y}) > c_1\;,\qquad
 \int_{\Gamma^K_{x,y}} \bigl(\DD{x}{y}^n(\omega)\bigr)^{-2}\,\PP^n(d\omega) < c_2\;,
\end{equ}
hold for every $n \ge 0$, whenever $\tilde V(x,y)\le K$. The integral over $\Gamma^K_{x,y}$ in \eref{e:assProb} should be interpreted as the integral over the projection of $\Gamma^K_{x,y}$ onto its $n$ first components.
\end{assumption}

A typical choice for $\Gamma^K_{x,y}$ is $\Gamma^K_{x,y} = A_{y,k}$ or $\Gamma^K_{x,y} = A_{x,k} \cap A_{y,k}$ with $k$ sufficiently large as a function of $K$. In this case, \lem{lem:bound} ensures that the conditions required on $\Gamma^K_{x,y}$ are satisfied.
As a consequence of Assumption~\ref{ass:Prob}, we have

\begin{proposition}\label{prop:coupl}
Let $\pQ^\infty_{x,y}$ be defined as above and suppose that assumptions~\ref{ass:Lyap} and \ref{ass:Prob} hold. Then there exists for every $K$ a constant $C_K$ such that
$\|\pQ^\infty_{x,y}\|\ge C_K$, whenever $\tilde V(x,y)\le K$.
\end{proposition}

\begin{proof}
Notice first that if $\mu_1$ and $\mu_2$ are two equivalent probability measures with $\mu_2(dx) = \Dens(x)\,\mu_1(dx)$, then the condition
\begin{equ}
 \int_{\!A} \bigl(\Dens(x)\bigr)^{-2}\,\mu_1(dx) < c
\end{equ}
implies that
\begin{equ}
 \bigl(\mu_1\wedge\mu_2\bigr)(A) \ge {\mu_1(A)^2\over 4c}\;,
\end{equ}
see, \eg \cite{MatNS}. Recalling \eref{e:relDens}, we use \lem{lem:bound} and the above estimate with $\mu_1 = \PP^n$, $\Dens = \Dens_{x,y}^n$, and $A = \Gamma^K_{x,y}$. Taking the limit $n\to\infty$ and using the assumption on $\Gamma^K_{x,y}$ proves the claim.
\end{proof}

Our last assumption will ensure that trajectories that have already coupled for some time have a very strong tendency to couple for all times.

In order to formulate our assumption, we introduce a family of sets $\sC^n_K(x,y)$, which are the possible final states of a ``coupled'' trajectory of length $n$, starting from $(x,y)$, and never leaving the set $\{(a,b)\;|\;V(a) + V(b) \le K\}$.
For a given pair of initial conditions $(x,y)\in\XX^2$ with $\tilde V(x,y) \le K$, we define the family of sets $\sC^n_K(x,y)\subset\XX\times\XX$ recursively in the following way:
\begin{equs}
\sC^0_K(x,y) &= \{(x,y)\}\;,\\
\sC^{n+1}_K(x,y) &= \bigcup_{(a,b)\in \sC^{n}_K(x,y)}
  \bigl\{\bigl(\Phi_{a,b}\circ\dPsi{a}{b}\bigr)(\omega)\;\big|\;\omega\in\SP\;\text{and}\;\dW{a}{b}(\omega)\le K\bigr\}\;.
\end{equs}
Notice that we would have obtained the same sets by reversing the directions of the arrows in the definition.

We also denote by $\Dens_{x,y}(\omega)$ the density of $\dmap{x}{y}^*\PP$ relative to $\PP$.

\begin{assumption}\label{ass:Coupl}
There exist positive constants $C_2$, $\gamma_2$ and $\zeta$, such that for every $K > 0$, every $(x_0,y_0)\in \XX^2$ with $\tilde V(x_0,y_0) \le K$, and every $(x,y) \in \sC_K^n(x_0,y_0)$, the estimate
\begin{equ}[e:ass4]
 \int_{\iW{x}{y}(\omega) \le K}\bigl(1-\Dens_{x,y}(\omega)\bigr)^2\,\PP(d\omega) \le C_2 e^{-\gamma_2 n}(1+K)^\zeta\;,
\end{equ}
holds for $n > \zeta\ln(1+K)/\gamma_2$.
\end{assumption}

This assumption means that if the process couples for a time $n$, the density $\Dens_{x,y}$ is close to $1$ on an increasingly large set, and therefore the probability of coupling for a longer time becomes increasingly large. This assumption is sufficient for the family of measures $(\pR\pQ^n)_{x,y}$ to have an exponential tail at large values of $n$. More precisely, we have
\begin{proposition}\label{prop:decExp}
Let assumptions~\ref{ass:Lyap}, \ref{ass:W} and \ref{ass:Coupl} hold. Then, there exists a positive constant $\gamma_3$ and, for every $K>0$, a constant $C_K$ such that
\begin{equ}[e:estRQn]
  \bigl\|(\pR\pQ^n)_{x,y}\bigr\| \le C_K e^{-\gamma_3 n}\;,
\end{equ}
holds for every $n > 0$, whenever $\tilde V(x,y)\le K$.
\end{proposition}

We first show the following elementary estimate (it is not optimal, but sufficient for our needs):
\begin{lemma}\label{lem:estProb}
Let $\mu_1,\mu_2 \in\Meas_1(\XX)$ be two equivalent probability measures with $\mu_2(dx) = \Dens(x)\,\mu_1(dx)$. Then the conditions
\begin{equ}
\mu_1(A) \ge 1-\eps_1\quad\text{and}\quad \qquad \int_{\!A} \bigl(1-\Dens(x)\bigr)^{2}\,\mu_1(dx) \le \eps_2\;,
\end{equ}
for some measurable set $A$ imply that
\begin{equ}
 \bigl(\mu_1\wedge\mu_2\bigr)(A) \ge 1 - \eps_1 - \eps_2^{1/2}\;.
\end{equ}
\end{lemma}
\begin{proof}
Define the set $E\subset \XX$ by
\begin{equ}
 E = A \cap \{x\in\XX\;|\;\Dens(x) \ge 1\}\;.
\end{equ}
We then have
\begin{equs}
 \bigl(\mu_1\wedge\mu_2\bigr)(A) &= \mu_1(E) + \int_{A\setminus E} \Dens(x)\, \mu_1(dx)\\
&= \mu_1(A) - \int_{A\setminus E} \bigl(1-\Dens(x)\bigr)\, \mu_1(dx)\\
&\ge \mu_1(A) - \int_{A\setminus E} \bigl|1-\Dens(x)\bigr|\, \mu_1(dx)\\
&\ge 1-\eps_1 - \sqrt{\int_{A\setminus E} \bigl(1-\Dens(x)\bigr)^2\, \mu_1(dx)}\;.
\end{equs}
This shows the claim.
\end{proof}

\begin{proof}[of \prop{prop:decExp}]
Fix the value $n$ and the pair $(x,y)$. For every $c_n \ge \tilde V(x,y)$ (we will fix it later), we have the estimate
\begin{equs}\label{e:estRQ}
 \bigl\|(\pR\pQ^n)_{x,y}\bigr\| &= \int_{\XX^2} \bigl(1-\|\dmap{x_n}{y_n}^*\PP\wedge\PP\|\bigr)\, \bigl(\pi_n^*\pQ^n_{x,y}\bigr)(dx_n,dy_n)\\
  &\le \bigl(\pi_n^*\pQ^n_{x,y}\bigr)\bigl(\XX^2\setminus\sC_{c_n}^n(x,y)\bigr) \\ 
	&\qquad +\int_{\sC_{c_n}^n(x,y)} \bigl(1 - \bigl\|\dmap{x_n}{y_n}^*\PP\wedge\PP\bigr\|\bigr)\, \bigl(\pi_n^*\pQ^n_{x,y}\bigr)(dx_n,dy_n)\;.
\end{equs}
Now choose another value $w_n$ to be fixed later and consider for every $(x_n,y_n)$ the set
\begin{equ}
B_n = \bigl\{\omega \in \SP\,|\, \iW{x_n}{y_n}(\omega) \le w_n\bigr\} \;.
\end{equ}
By the definition of $\sC_{c_n}^n(x,y)$, its elements $(x_n,y_n)$ satisfy in particular the inequality $\tilde V(x_n,y_n) \le c_n$.
By Assumption~\ref{ass:W} and the Lyapunov structure, we have for every $(x_n,y_n) \in \sC_{c_n}^n(x,y)$ the estimate
\begin{equ}
 \PP(B_n) \ge 1 - C{c_n \over w_n^{1/\delta}}\;.
\end{equ}
Combining this and Assumption~\ref{ass:Coupl} with \lem{lem:estProb} yields
\begin{equ}
 \bigl\|\dmap{x_n}{y_n}^*\PP\wedge\PP\bigr\| \ge 1 - C{c_n \over w_n^{1/\delta}} - C e^{-\gamma_2 n /2} \bigl(1 + w_n\bigr)^{\zeta/2}\;,
\end{equ}
as long as $w_n$ is such that
\begin{equ}[e:condn]
w_n \ge c_n \quad
\text{and}\quad n \ge \zeta \ln(1 + w_n)/\gamma_2\;.
\end{equ}
It remains to give an upper bound for $\bigl(\pi_n^*\pQ^n_{x,y}\bigr)\bigl(\XX^2\setminus\sC_{c_n}^n(x,y)\bigr)$ to complete our argument. Define the sets $A^n(K) \subset \XX^n\times\XX^n\times\SP^n\times\SP^n$ by
\begin{equ}
 A^n(K) = \bigl\{(x_i,y_i,\omega_i,\eta_i)_{i=1}^{n} \;|\; W(x_i,\omega_i) + W(y_i,\eta_i)\le K\bigr\}\;.
\end{equ}
It is clear by the definition of $\sC_{c_n}^n(x,y)$ that we have the equality
\begin{equ}
 \bigl(\pi_n^*\pQ^n_{x,y}\bigr)\bigl(\XX^2\setminus\sC_{c_n}^n(x,y)\bigr) = \pQ^n_{x,y}\bigl(\XX^n\times\XX^n\times\SP^n\times\SP^n\setminus A^n(c_n)\bigr)\;,
\end{equ}
where $\pQ^n_{x,y}$ is considered as a measure on $\XX^n\times\XX^n\times\SP^n\times\SP^n$, following \rem{rem:prob}. Since $\pQ^n_{x,y} \le \pC^n_{x,y}$, we have
\begin{equ}
 \bigl(\pi_n^*\pQ^n_{x,y}\bigr)\bigl(\XX^2\setminus\sC_{c_n}^n(x,y)\bigr) \le 1 - \pC^n_{x,y}\bigl(A^n(c_n)\bigr) \le C{n \bigl(\tilde V(x,y) + 1\bigr) \over c_n}\;,
\end{equ}
for some constant $C$. This last estimate is obtained in a straightforward way, following the lines of the proof of \lem{lem:bound}.
Plugging these estimates back into \eref{e:estRQ} yields
\begin{equ}
 \bigl\|(\pR\pQ^n)_{x,y}\bigr\| \le C{n \bigl(\tilde V(x,y) + 1\bigr) \over c_n}  + C{c_n \over w_n^{1/\delta}} + C e^{-\gamma_2 n /2} \bigl(1 + w_n\bigr)^{\zeta/2}\;.
\end{equ}
At this point, we make use of our freedom to choose $c_n$ and $w_n$. We set
\begin{equs}
 c_n = \tilde V(x,y) + e^{\gamma_c n}\qquad\text{and}\qquad
w_n =   \tilde V(x,y) + e^{\gamma_w n}\;,
\end{equs}
with $\gamma_c$ and $\gamma_w$ given by
\begin{equs}
 \gamma_c = {1\over 2+2\delta\zeta}\gamma_2\qquad\text{and}\qquad
\gamma_w = {\delta\over 1+\delta\zeta}\gamma_2\;.
\end{equs}
As a consequence,
there exist for any $\gamma < \gamma_c$ some constants $C$ and $c$ such that
\begin{equ}
 \bigl\|(\pR\pQ^n)_{x,y}\bigr\| \le C \bigl(1 + \tilde V(x,y)\bigr)^c e^{-\gamma n}\;,
\end{equ}
as long as $n \ge \zeta\ln(1+w_n)/\gamma_2$. (Such a value of $n$ can always be found, because the exponent $\gamma_w$ is always smaller than $\gamma_2/\zeta$.)
In order to complete the argument, we notice that \eref{e:estRQn} is trivially satisfied for small values of $n$ because $\bigl\|(\pR\pQ^n)_{x,y}\bigr\|$ is always smaller than $1$ by definition: it suffices to choose $C_K$ sufficiently big.
The proof of \prop{prop:decExp} is complete.
\end{proof}

\section{An Exponential Mixing Result}
\label{sect:abstr}

This section is devoted to the proof of the main theorem of this paper. 

\begin{theorem}\label{theo:main}
Let $\PH$ be a RDS with state space $\XX$ satisfying assumptions~\ref{ass:Lyap}--\ref{ass:Coupl}. Then, there exists a constant $\gamma>0$ such that
\begin{equ}
 \bigl\|\pP^n_x - \pP^n_y\|_\Lip \le C\bigl(1+\tilde V(x,y)\bigr)\,e^{-\gamma n}\;,
\end{equ}
for every $(x,y) \in \XX^2$ and every $n > 0$.
\end{theorem}

\begin{remark}
The proof of \theo{theo:main} does not rely on assumptions~\ref{ass:Prob} and \ref{ass:Coupl} directly, but on the conclusions of Propositions~\ref{prop:coupl} and \ref{prop:decExp}. Nevertheless, in the setting of stochastic differential equations, it seems to be easier to verify the assumptions rather than to show the conclusions of the propositions by other means.
\end{remark}

\begin{corollary}\label{cor:main}
If $\PH$ satisfies assumptions~\ref{ass:Lyap}--\ref{ass:Coupl}, it possesses a unique invariant measure $\mu_*$ and
\begin{equ}
 \bigl\|\pP^n_x - \mu_*\|_\Lip \le C\bigl(1+V(x)\bigr)\, e^{-\gamma n}\;.
\end{equ}
\end{corollary}

\begin{proof}[of the corollary]
To show the existence of the invariant measure $\mu_*$, we show that for any given initial condition $x$ with $V(x)<\infty$, the sequence of measures $\pP^n_x$ is a Cauchy sequence in the norm $\|\cdot\|_\Lip$. We have indeed
\begin{equs}
 \|\pP^n_x -\pP^{n+k}_x\|_\Lip &= \sup_{\|g\|_\Lip \le 1} \int_\XX g(z)\bigl(\pP^n_x -\pP^{n+k}_x\bigr)(dz)\\
&= \sup_{\|g\|_\Lip \le 1} \int_\XX \int_\XX g(z)\bigl(\pP^n_x -\pP^{n}_y\bigr)(dz)\,\pP^k_x(dy)\\
&\le \int_\XX \|\pP^n_x -\pP^{n}_y\|_\Lip\,\pP^{k}_x(dy)
\le C e^{-\gamma n}\int_\XX \bigl(1+\tilde V(x,y)\bigr)\,\pP^k_x(dy) \\
&\le C e^{-\gamma n} \bigl(1+V(x)\bigr)\;,
\end{equs} 
where we used the Lyapunov structure to get the last inequality.

The claim now follows immediately from the theorem, noticing that if $\mu_*$ is an invariant measure for $\PH$, then
\begin{equ}
 \int_\XX V(x)\,\mu_*(dx) \le {b \over 1-a}\;,
\end{equ}
due to the Lyapunov structure and the fact that the dynamics immediately leaves the set $V^{-1}(\infty)$.
\end{proof}

Before we turn to the proof of \theo{theo:main}, we introduce some notations and make a few remarks. By iterating \eref{e:defCxy}, one sees that
\begin{equ}[e:equiv]
 \pC^\infty_{x,y} = \pQ^\infty_{x,y} + \sum_{n=0}^\infty \bigl(\pC^\infty \pR \pQ^n\bigr)_{x,y}\;,
\end{equ}
where the symbol $\bigl(\pC^\infty \pR \pQ^n\bigr)_{x,y}$ is to be interpreted in the sense of \eref{e:comp}. This expression is the equivalent, in our setting, of Lemma~2.1 in \cite{MatNS}. Using \eref{e:equiv}, the Markov chain $\Psi$ can be described by means of another Markov chain $\Upsilon$ on $\YY = (\XX^2 \times \N)\cup\{\star\}$, where $\star$ corresponds to ``coupling for all times'' in the sense of \sect{sect:couplconstr}.
First, we define
\begin{equ}[e:defK0]
 K_0 = {4b \over 1-a}\;,\qquad \tilde K_0 = \bigl\{(x,y)\;|\; \tilde V(x,y) \le K_0\bigr\}\;,
\end{equ}
where $a$ and $b$ are the constants appearing in the Lyapunov condition. This set is chosen in such a way that 
\begin{equ}[e:propK0]
\int_{\XX\times\XX} \tilde V(x,y) \,\pC_{x_0,y_0}(dx,dy) \le {1 + a \over 2}\tilde V(x_0,y_0)\;,\qquad \forall\; (x_0,y_0) \not\in \tilde K_0\;.
\end{equ}

At time $0$, $\Upsilon$ is located at $(x,y,0)$.
If it is located at $(x,y,n)$ and $(x,y) \not\in \tilde K_0$, then it makes one step according to $\pC_{x,y}$ and $n$ is incremented by one:
\begin{equ}
  \nP^_{(x,y,n)} = \pC_{x,y} \times \delta_{n+1}\;.
\end{equ}
If $\Upsilon$ is located at $(x,y,n)$ and $(x,y) \in \tilde K_0$, then it has a probability $\|\pQ_{x,y}^\infty\|$ of jumping to $\star$ and a probability $\|(\pR\pQ^m)_{x,y}\|$ of making $m$ steps according to $(\pR\pQ^m)_{x,y}$:
\begin{equ}
  \nP^_{(x,y,n)} = \|\pQ_{x,y}^\infty\|\delta_\star + \sum_{m=0}^\infty \pi_m^*(\pR\pQ^m)_{x,y}\times\delta_{n+m+1}\;.
\end{equ}
If $\Upsilon$ is located at $\star$, it remains there:
\begin{equ}
 \nP^_\star = \delta_\star \;.
\end{equ}
The Markov chain $\Upsilon$ induces a family $\hPt^\infty_{x,y}$ of probability measures on $\YY^\infty$.
Let $\hat\tau\from\YY^\infty \to \N\cup\{\infty\}$ be the function that associates to a sequence of elements in $\YY$ the largest value of $n$ that is reached by the sequence ($\hat\tau = 0$ if the sequence is equal to $\star$ repeated). We also define $\hat\kappa\from\YY^\infty \to \N\cup\{\infty\}$ as the value of $n$ attained at the first non-vanishing time when the sequence hits the set $\tilde K_0 \times \N$ ($\hat\kappa = \infty$ if this set is never reached). The construction of $\Upsilon$ is very close to the coupling construction of \cite{MatNS}.

The crucial observation for the proof of \theo{theo:main} is
\begin{lemma}
Let $\PH$ be a RDS with state space $\XX$ satisfying assumptions~\ref{ass:Lyap} and \ref{ass:Conv}, and let $\Upsilon$ be defined as above. Then, there exists a constant $C$ such that
\begin{equ}
 \bigl\|\pP^n_x - \pP^n_y\|_\Lip \le \hPt^\infty_{x,y}\bigl(\{\hat\tau \ge n/2\}\bigr) + C e^{-\gamma_1 n/2}\;,
\end{equ}
for every $(x,y) \in \XX^2$ and every $n > 0$.
\end{lemma}

\begin{proof}
Recall the Markov chain $\hat\Psi$ defined in \rem{rem:Psihat}. We define
a function $\tau_1$ on its pathspace by
\begin{equs}
\tau_1\from\XX^\infty\times\XX^\infty\times\{0,1\}^\infty &\to \N\cup\{\infty\} \\
 \bigl\{(x_i,y_i,b_i)\bigr\}_{i=1}^\infty &\mapsto \inf \bigl\{n\,|\,(x_n,y_n)\in\tilde K_0\;\text{and}\;b_i = 1\;\forall i \ge n\bigr\}\;.
\end{equs}
Combining \eref{e:couplingDiff} with Assumption~\ref{ass:Coupl} and the
definition of $\tau_1$, one sees that
\begin{equ}
  \bigl\|\pP^n_x - \pP^n_y\|_\Lip \le \pC^\infty_{x,y}\bigl(\{\tau_1 \ge n/2\}\bigr) + C e^{-\gamma_1 n/2}\;.
\end{equ}
From the construction of $\Upsilon$ and the definition of $\hat\tau$, we see furthermore that the probability distributions of $\tau_1$ under $\pC^\infty_{x,y}$ and of $\hat\tau$ under $\hPt^\infty_{x,y}$ are the same.
\end{proof}

\begin{proof}[of \theo{theo:main}]
It remains to show that $\hPt^\infty_{x,y}\bigl(\{\hat\tau \ge n/2\}\bigr)$ has an exponential tail. The key observation is the following. Let $x_n \in \N\cup\{-\infty\}$ with $n\ge 0$ be a Markov chain defined by
\begin{equ}
x_0 = 0\;,\quad x_{n+1} = \cases{-\infty& with probability $p_\star$,\cr
	x_n + m& with probability $p_m$,}
\end{equ}
where $m\ge 1$ and, of course, $p_\star + \sum_{m=1}^\infty p_m = 1$. 
\begin{lemma}\label{lem:exp}
If the $p_m$ have an exponential tail and we define $\tau = \max_n x_n$, then the probability distribution of $\tau$ also has an exponential tail. 
\end{lemma}
\begin{proof}
The claim is an easy consequence of Kendall's theorem, but for the sake of completeness, and because the proof is quite elegant, we outline it here.
Define the analytic function $p(\zeta) = \sum_{m=1}^\infty p_m \zeta^m$
and define $q_n$ as the probability of $\tau$ being equal to $n$. Notice that,
because of the exponential tail, $p$ is analytic in a disk of radius strictly bigger than $1$ around the origin.
A straightforward computation shows that $q_0 = p_\star$ and, for $n > 0$,
\begin{equ}
 q_n = p_\star \Bigl(p_n +\sum_{k_1 + k_2 = n}p_{k_1}p_{k_2}
	+\sum_{k_1 + k_2 + k_3 = n}p_{k_1}p_{k_2}p_{k_3} + \ldots\Bigr)\;,
\end{equ}
which is equal to the $n$th Taylor coefficient of the function
\begin{equ}
 q(\zeta) = {p_\star \over 1 - p(\zeta)}\;.
\end{equ}
Since $p(1) = 1-p_\star < 1$, there exists an $\eps > 0$ such that
$p(1 + \eps) < 1$. Furthermore, since the $p_n$ are all positive, one has the estimate $|p(\zeta)| \le p(|\zeta|)$. Using Cauchy's formula on a circle of radius $1+\eps$, one gets
\begin{equ}
 |q_n| \le {p_\star\over 1-p(1+\eps)}{1\over (1+\eps)^n}\;,
\end{equ}
which shows the claim.
\end{proof}
Before we prove \theo{theo:main} in full generality, we restrict ourselves
to the case when $(x,y) \in \tilde K_0$. It follows from the construction that
$\hat \tau$ (seen as a random variable under the distribution induced by $\hPt^\infty_{x,y}$) is dominated by the process $x_n$ constructed above with the tail distribution of the $p_m$ being equal to
\begin{equ}
 \tilde p_m = \sup_{(x,y) \in \tilde K_0} \hPt^\infty_{x,y}\bigl(\bigl\{\hat \kappa = m\bigr\}\bigr)\;.
\end{equ}
This means that we define $m_*$ as
\begin{equ}
 m_* = \inf \Bigl\{m\;\Big|\;\sum_{n=m}^\infty \tilde p_n \le 1\Bigr\}\;,
\end{equ}
and then set $p_m = \tilde p_m$ for $m \ge m_*$, $p_m = 0$ for $m < m_*-1$,
and $p_{m_*-1}$ in such a way that the $p_n$ sum up to $1$.

Because of \lem{lem:exp}, it suffices to show that the tail distribution of the $\tilde p_m$ decays exponentially. We thus estimate the quantity
$\hPt^\infty_{x,y}\bigl(\bigl\{\hat \kappa \ge n\bigr\}\bigr)$. To this end, we introduce the function $\tau_\Psi^{}\from\XX^\infty\times\XX^\infty\to\N\cup\{\infty\}$ defined by
\begin{equ}
 \tau_\Psi^{} (x,y) = \inf\{n > 0\;|\; (x_n,y_n)\in\tilde K_0\}\;.
\end{equ}

Notice that, in order to have $\hat\kappa \ge n$, there are two possibilities. Either the first step of $\Upsilon$ is taken according to $(\pR\pQ^m)_{x,y}$ with some $m\ge n/2$, or the corresponding realization of $\Psi$ stays outside of
$\tilde K_0$ for a time longer than $n/2$. This yields the estimate
\begin{equ}
 \hPt^\infty_{x,y}\bigl(\bigl\{\hat \kappa \ge n\bigr\}\bigr)
\le \sum_{m=n/2}^\infty \bigl\|(\pR\pQ^m)_{x,y}\bigr\| + {n\over 2}\sup_{(x_0,y_0)\in\tilde K_0} \pC_{x_0,y_0}^\infty\bigl(\bigl\{\tau_\Psi^{} \ge n/2\bigr\}\bigr)\;,
\end{equ}
holding for $(x,y)\in\tilde K_0$.
The first term has an exponential tail by \prop{prop:decExp}. The second term has also an exponential tail by \eref{e:propK0} and standard Lyapunov techniques (see \eg \cite[Thm~15.2.5]{MT}). This concludes the proof of \theo{theo:main} for the case $(x,y)\in \tilde K_0$.

In order to conclude the proof for the case $(x,y)\not\in \tilde K_0$, notice that
\begin{equs}
 \hPt^\infty_{x,y}\bigl(\bigl\{\hat\tau \ge n\bigr\}\bigr) &\le \sum_{m = 1}^\infty \pC_{x,y}^\infty\bigl(\bigl\{\tau_\Psi = m\bigr\}\bigr) \sup_{(x_0,y_0)\in\tilde K_0} \hPt^\infty_{x_0,y_0}\bigl(\bigl\{\hat\tau \ge n-m\bigr\}\bigr) \\
&\le {n \over 2}\sup_{(x_0,y_0)\in\tilde K_0} \hPt^\infty_{x_0,y_0}\bigl(\bigl\{\hat\tau \ge n/2\bigr\}\bigr) + \sum_{m = n/2}^\infty \pC_{x,y}^\infty\bigl(\bigl\{\tau_\Psi = m\bigr\}\bigr)\;.
\end{equs}
The first term is bounded by the construction above.
The Lyapunov structure implies that there exists a constant $\gamma>0$ such that the first hitting time $\tau_\Psi$ satisfies $\Exp_{(x,y)} e^{\gamma \tau_\Psi} = \CO(\tilde V(x,y))$ for every $(x,y)\in\XX^2$ (see again \cite[Thm.~15.2.5]{MT}). This allows to bound the second term and concludes the proof of \theo{theo:main}.
\end{proof}

\section{Application to Stochastic Differential Equations}
\label{sect:equadiff}

In this section, we will see how to apply \theo{theo:main} to the case
when the RDS $\PH$ is constructed by sampling the solution of a (possibly
infinite-dimensional) stochastic differential equation. We will restrict ourselves to the case where the equation is driven by additive white noise. The case of multiplicative noise requires further estimates, but can also be described by the formalism exposed here.

Consider the equation described by
\begin{equ}[e:equaDiff]
 dx(t) = Ax\,dt + F(x)\,dt + Q\,d\omega(t)\;,\qquad x(0)=x_0\;,
\end{equ}
where $x$ belongs to some separable Hilbert space $\CH$, $\omega$ is the cylindrical Wiener process on some separable Hilbert space $\CW$, and $A$, $F$ and $Q$ satisfy the following assumptions:
\begin{bassumption}\label{ass:ExistSol}
\begin{claim}
\item[{\it a.}] The linear operator $A\from\Dom(A)\to\CH$ is the generator of a strongly continuous semigroup on $\CH$.
\item[{\it b.}] The operator $e^{At}Q\from\CW\to\CH$ is Hilbert-Schmidt for every $t>0$ and one has the estimate
\begin{equ}[e:condQ]
 \int_0^1 \bigl\|e^{At}Q\bigr\|^2_\HS\,dt <\infty\;.
\end{equ}
\item[{\it c.}] The nonlinear operator $F\from\Dom(F)\to\CH$ is such that, for every $x_0\in\CH$, there exists a unique, continuous stochastic process $x(t)$ such that $x(s) \in \Dom(F)$ for $s > 0$ and
\begin{equ}[e:defSol]
 x(t) = e^{At} x_0 + \int_0^t e^{A(t-s)}F\bigl(x(s)\bigr)\,ds + \int_0^t e^{A(t-s)}Q\,d\omega(s)\;,
\end{equ}
for every $t>0$.
\end{claim}
\end{bassumption}
\begin{remark}
This assumptions simply states that there exists a unique weak solution to \eref{e:equaDiff}. Notice that we do {\em not} make any assumptions on the tightness of the transition probabilities for \eref{e:equaDiff}. As a consequence, existence and uniqueness results for invariant measures can in principle be deduced from \theo{theo:mainEqu} below even in cases where the semigroup $e^{At}$ is not compact.
\end{remark}

In order to recover the formalism used in \sect{sect:coupl}, we follow \cite{ZDP1} and introduce an auxiliary Hilbert space $\hat\CW$ such that there exists a continuous embedding $\iota\from\CW \hookrightarrow\hat\CW$, which is Hilbert-Schmidt. We can now set $\SP = \CC_0([0,1],\hat\CW)$, the space of continuous $\hat\CW$-valued functions that vanish at $0$, and define $\PP$ as the Wiener measure on $\hat\CW$ with covariance operator $\iota\iota^*$.

We define $\PH\from\CH\times\SP \to \CH$ as the map that solves \eref{e:defSol} up to time $1$ given an initial condition and a realization of the noise. This map is defined $\PP$-almost everywhere on $\SP$. We also denote by $\Phi_t\from\CH\times\SP^\infty\to\CH$ the map that maps an initial condition and a realization of the noise onto the solution of \eref{e:defSol} after a time $t$.

Our next assumption is the existence of an appropriate Lyapunov function $V$:
\begin{bassumption}\label{ass:Lyap2}
There exists a measurable function $V\from\CH\to[0,\infty]$ and constants $a < 1$ and $b,c,d>0$ such that
\begin{equs}
\Exp_\omega V\bigl(\PH(x,\omega)\bigr) &\le aV(x) + b\;,\\
\Exp_\omega \Bigl(\sup_{0 \le t \le 1}V\bigl(\Phi_t(x,\omega)\bigr) \Bigr)&\le cV(x) + d\;,\label{e:Lyap2}\\
\PP\bigl(\bigl\{\omega\;|\;V\bigl(\PH(x,\omega)\bigr) = \infty\bigr\}\bigr) &= 0\;,
\end{equs}
for every $x\in\CH$. Furthermore, $V$ dominates the norm in $\CH$ in the sense that $\|x\| \le C\bigl(1+V(x)\bigr)$ for some constant $C$.
\end{bassumption}
As is \sect{sect:Prop}, we define $\tilde V(x,y) = V(x) + V(y)$.
\begin{remark}
Take $\CH$ equal to $\L^2(\CO)$ for some regular bounded domain $\CO \subset \R^d$, $A$ a second-order elliptic differential operator on $\CO$ with sufficiently smooth coefficients, and $F$ any polynomial non-linearity of odd degree having the correct sign. The assumptions~\ref{ass:ExistSol} and \ref{ass:Lyap2} are satisfied with $V(x) = \|x\|_\dstar^p$ for every power $p\ge 1$ and every ``reasonable'' norm $\|\cdot\|_\dstar^{}$, as long as $Q$ is ``small'' enough. (One can for example take for $\|\cdot\|_\dstar^{}$ the $\L^\infty$ norm or a Sobolev norm.)
\end{remark}

We now turn to the binding construction for the problem \eref{e:equaDiff}.
Take a function $G\from\CH^2 \to \CW$ and consider the $\CH^2$-valued process $(x,y)$ solving
\minilab{e:defcoupEqu}
\begin{equs}
	dx(t) &= Ax\,dt + F(x)\,dt + Q\,d\omega(t)\;,\label{e:defx}\\
  dy(t) &= Ay\,dt + F(y)\,dt + Q\,G(x,y)\,dt + Q\,d\omega(t)\;.\label{e:defy}
\end{equs}
Notice that the realization of $\omega$ is the same for both components.
The process \eref{e:defcoupEqu} yields our binding construction for \eref{e:equaDiff}.
In order to give sense to \eref{e:defy}, we introduce the $\CH$-valued process $\rho(t) = y(t) - x(t)$ and we define it pathwise as the solution of the ordinary differential equation
\begin{equ}[e:defRho]
 \dot \rho = A\rho + F(x + \rho) - F(x) + Q\,G(x,x +\rho)\;.
\end{equ}
We assume that $G$ is sufficiently regular to ensure the existence and uniqueness of global weak solutions to \eref{e:defRho} for almost every (with respect to the measure on pathspace induced by $\PH_t$) continuous function
$x \from [0,\infty) \to \CH$. This allows us to {\it define} the stochastic process $y(t)$ by $y(t) = x(t) + \rho(t)$. We will denote by $\vec\Phi_t\from\XX\times\XX\times\SP\to\XX$ the map that solves \eref{e:defy} up to time $t$, given an initial condition for $x$ and $y$, and a realization of the noise.

The above construction is invertible in the following sense. Consider the $\CH^2$-valued process
\minilab{e:defcoupEquInv}
\begin{equs}
	d\tilde x(t) &= A\tilde x\,dt + F(\tilde x)\,dt - Q\,G(\tilde x,\tilde y)\,dt  + Q\,d\tilde\omega(t)\;,\label{e:defxtilde}\\
  d\tilde y(t) &= A\tilde y\,dt + F(\tilde y)\,dt + Q\,d\tilde\omega(t)\;,\label{e:defytilde}
\end{equs}
where we give sense to the equation for $\tilde x$ as above by setting $\tilde \rho = \tilde y - \tilde x$ and solving
\begin{equ}
  \dot {\kern-.2em{\tilde \rho}} = A\tilde \rho + F(\tilde y) - F(\tilde y - \tilde\rho) + Q\,G(\tilde y - \tilde \rho,\tilde y)\;.
\end{equ}
We denote by $\cev\Phi_t\from\XX\times\XX\times\SP\to\XX$ the map that solves \eref{e:defxtilde} up to time $t$, given an initial condition for $\tilde x$ and $\tilde y$, and a realization of the noise $\tilde\omega\in\Omega$. We see that \eref{e:defcoupEquInv} can be obtained from \eref{e:defcoupEqu} by the substitution $d\tilde\omega = d\omega + G(x,y)\,dt$ and a renaming of the variables. This observation yields the invertibility of the maps $\dmap{x}{y}$ defined in Eq.~\eref{e:bindEqu} below.

We will state two more assumptions to make sure that the conclusions of \theo{theo:main} hold. 
First, we want $G$ to become small as $x$ and $y$ become close.

\begin{bassumption}\label{ass:ExpSmall}
There exists a constant $C>0$ and exponents $\alpha,\beta>0$ such that
\begin{equ}[e:boundG]
 \|G(x,y)\|^2 \le C \|x - y\|^\alpha \bigl(1 + \tilde V(x,y)\bigr)^{\beta}\;,
\end{equ}
for every $x,y \in \CH$.
\end{bassumption}
The last assumption ensures that the process $y(t)$ converges towards $x(t)$ for large times. 
\begin{bassumption}\label{ass:Conv2}
There exist positive constants $C$ and $\gamma$ such that the solutions of \eref{e:defcoupEqu} and \eref{e:defcoupEquInv} satisfy
\minilab{e:boundrho}
\begin{equs}
 \|\Phi_t(x,\omega) - \vec\Phi_t(x,y,\omega)\| &\le C e^{-\gamma t}\Bigl(1 +V(y) + \sup_{s \le t} V\bigl(\Phi_s(x,\omega)\bigr) \Bigr)\;,\quad\qquad\label{e:boundrhodir}\\
 \|\cev\Phi_t(x,y,\omega) - \Phi_t(y,\omega)\| &\le C e^{-\gamma t}\Bigl(1 + V(x) + \sup_{s \le t} V\bigl(\Phi_s(y,\omega)\bigr) \Bigr)\;,\quad\qquad\label{e:boundrhoinv}
\end{equs}
for $\PP^\infty$-almost every $\omega\in\SP^\infty$. Furthermore, there exists $\delta > 0$ such that one has the estimate%
\minilab{e:boundmodif}
\begin{equs}
 V\bigl(\vec\Phi_t(x,y,\omega)\bigr) &\le  C\Bigl(1+ V(y) + \sup_{s \le t} V\bigl(\Phi_s(x,\omega)\bigr)\Bigr)^\delta\;,\qquad\label{e:boundmodif1}\\
 V\bigl(\cev\Phi_t(x,y,\omega)\bigr) &\le  C\Bigl(1+ V(x) +  \sup_{s \le t} V\bigl(\Phi_s(y,\omega)\bigr)\Bigr)^\delta\;,\qquad\label{e:boundmodif2}
\end{equs}
for $\PP^\infty$-almost every $\omega\in\SP^\infty$ and every $t\ge 0$.
\end{bassumption}

\begin{remark}\label{rem:boundV}
One important particular case is the choice $V(x) = \|x\|^p$, where the power $p$ is chosen in such a way that \eref{e:boundrho} is satisfied. Notice that in this case, the estimates \eref{e:boundmodif} are a straightforward consequence of \eref{e:boundrho}.

The function $G$ is then only required to satisfy a bound of the type
\begin{equ}
 \|G(x,y)\|^2 \le C\|x-y\|^\alpha\bigl(1 + \|x\| + \|y\|\bigr)^q\;,
\end{equ}
with $\alpha$ and $q$ some arbitrary positive exponents.

It is also possible to choose $V(x) = \|x\|_\star^p$, with $\|\cdot\|_\star$ the norm of some Banach space $\CB\subset\CH$. In this case, \eref{e:boundrho} with the $\CB$-norm replacing the $\CH$-norm in the left-hand side implies \eref{e:boundmodif}.
\end{remark}

\begin{remark}
An equivalent way of writing \eref{e:boundrhoinv} is
\begin{equ}[e:alternative]
 \|\Phi_t(x,\omega) - \vec\Phi_t(x,y,\omega)\| \le C e^{-\gamma t}\Bigl(1 +V(x) + \sup_{s \le t} V\bigl(\vec\Phi_s(x,y,\omega)\bigr) \Bigr)\;.
\end{equ}
The equation \eref{e:alternative} will be more natural in our examples, but \eref{e:boundrho} is more symmetric and more convenient for the proof of \theo{theo:mainEqu} below.
\end{remark}

All these assumptions together ensure that exponential mixing takes place:

\begin{theorem}\label{theo:mainEqu}
Let $A$, $F$ and $Q$ be such that assumptions~\ref{ass:ExistSol} and \ref{ass:Lyap2} are satisfied. If there exists a function $G\from\CH^2\to\CW$ such that assumptions~\ref{ass:ExpSmall} and \ref{ass:Conv2} hold, then the solution of \eref{e:equaDiff} possesses a unique invariant measure $\mu_*$ and there exist constants $C,\gamma > 0$ such that
\begin{equ}
 \|\nP^n_x - \mu_*\|_\Lip \le C e^{-\gamma n}\bigl(1 + V(x)\bigr)\;.
\end{equ}
\end{theorem}

\begin{proof}
It suffices to show that assumptions~\ref{ass:Lyap}--\ref{ass:Coupl} hold.
Assumption~\ref{ass:Lyap} follows immediately from Assumption~\ref{ass:Lyap2}. In order to check the other assumptions, we define the various objects appearing in the previous sections. We have already seen that $\XX = \CH$, $\SP = \CC_0([0,1],\hat\CW)$, and $\PH$ is the solution of \eref{e:equaDiff} at time $1$.

We define the function $W\from\CH\times\SP\to[1,\infty]$ by
\begin{equ}
 W(x,\omega) = \sup_{t\in[0,1]} V\bigl(\PH_t(x,\omega)\bigr)\;.
\end{equ}
The estimate \eref{e:Lyap2} and the definition ensure that $W$ satisfies \eref{e:propW1} and \eref{e:propW2}. The bound \eref{e:boundmodif} ensures that Assumption~\ref{ass:W} is also satisfied.

It remains to define the binding functions $\dmap{x}{y}$ and to compute the densities $\Dens^n_{x,y}$. According to the constructions \eref{e:defcoupEqu} and \eref{e:defcoupEquInv}, we define for $(x,y)\in\CH^2$ the binding functions
\minilab{e:bindEqu}
\begin{equs}
 \bigl(\dmap{x}{y}(\omega)\bigr)(t) &= \omega(t) + \int_0^t G\bigl(\Phi_s(x,\omega),\vec\Phi_s(x,y,\omega)\bigr)\,ds\;,\label{e:bindDir}\\
 \bigl(\imap{x}{y}(\omega)\bigr)(t) &= \omega(t) - \int_0^t G\bigl(\cev\Phi_s(x,y,\omega),\Phi_s(y,\omega)\bigr)\,ds\;,\label{e:bindInv}
\end{equs}
with $t\in[0,1]$.
It follows from the construction that these functions are each other's inverse. Furthermore, if we identify $\SP^n$ with $\CC_0([0,n],\hat\CW)$ in a natural way, we see that the maps $\Xi^n_{x,y}$ introduced in \eref{e:defXi} are obtained from \eref{e:bindEqu} by simply letting $t$ take values in $[0,n]$. These observations allow us to compute the densities $\Dens^n_{x,y}(\omega)$ by Girsanov's theorem:
\begin{lemma}\label{lem:Dens}
The family of densities $\Dens^n_{x,y}(\omega)$ is given by
\begin{equ}
 \Dens^n_{x,y}(\omega) = \exp\Bigl(\int_0^n G\bigl(\cev\Phi_t(x,y,\omega),\Phi_t(y,\omega)\bigr)\,d\omega(t) - {1\over 2}\int_0^n \|G(\ldots)\|^2\,dt\Bigr)\;,
\end{equ}
where the arguments of $G$ in the second term are the same as in the first term.
\end{lemma}
\begin{proof}
If we can show that Girsanov's theorem applies to our situation,
 then it precisely states that
\begin{equ}
 \imap{x}{y}^*\hat\PP^n = \PP^n\;,
\end{equ}
with $\hat\PP^n(d\omega) = \Dens^n_{x,y}(\omega)\,\PP^n(d\omega)$, and $\Dens^n_{x,y}(\omega)$ defined as above. Applying $\dmap{x}{y}^*$ to both sides of the equality shows the result.

We now show that Girsanov's theorem can indeed be applied. By \cite[Thm~10.14]{ZDP1}, it suffices to verify that
\begin{equ}[e:normD]
 \int_{\SP^n}\Dens^n_{x,y}(\omega)\,\PP^n(d\omega) = 1\;.
\end{equ}
This can be achieved by a suitable cut-off procedure. Define for $N>0$ the function
\begin{equ}
 \GN(x,y) = \cases{G(x,y)&if $\|G(x,y)\| \le N$,\cr 0&otherwise,}
\end{equ}
and define
\begin{equ}
  \Dens^{n,N}_{x,y}(\omega) = \exp\Bigl(\int_0^n \GN\bigl(\cev\Phi_t(x,y,\omega),\Phi_t(y,\omega)\bigr)\,d\omega(t) - {1\over 2}\int_0^n \|\GN(\ldots)\|^2\,dt\Bigr)\;.
\end{equ}
It is immediate that \eref{e:normD} holds for $\Dens^{n,N}_{x,y}$. Furthermore, it follows from Assumption~\ref{ass:Conv2} that there exists a constant $C_{N}$ such that $\Dens^{n,N}_{x,y}(\omega) = \Dens^{n}_{x,y}(\omega)$ on the set
\begin{equ}
 \Gamma_{\!N} = \bigl\{\omega\in\CP^n \;|\; \tilde V\bigl(\cev\Phi_s(x,y,\omega),\Phi_s(y,\omega)\bigr) < C_{N}\quad\forall s \in [0,n] \bigr\}\;.
\end{equ}
The sets $\Gamma_{\!N}$ satisfy $\lim_{N\to\infty}\PP^n(\Gamma_{\!N}) = 1$ by \eref{e:Lyap2} and \eref{e:boundmodif2}. This shows that \eref{e:normD} holds.
Notice that the {\it a-priori} bounds of Assumption~\ref{ass:Conv2} were crucial in this step in order to apply Girsanov's theorem. The bound \eref{e:boundG} alone could lead to exploding solutions for which Girsanov's theorem does not apply.
\end{proof}
It is immediate that Assumption~\ref{ass:Conv} follows from Assumption~\ref{ass:Conv2} and the definition of the norm $\|\cdot\|_\Lip$.

We now turn to the verification of Assumption~\ref{ass:Prob}. Recalling the definition \eref{e:defAxk}, we see that in our case
\begin{equ}
 A_{y,k} \subset B_{y,k} \equiv \bigl\{\omega\in\SP^\infty \;|\; V\bigl(\Phi_s(y,\omega)\bigr) \le k\bigl(V(y) + s^2\bigr) \quad\forall s \ge 0 \bigr\}\;.
\end{equ}
As we see from the definition of $B_{y,k}$, a natural definition for a truncation $G_{y,k}$ of $G$ (this time the truncation additionally depends on time) is 
\begin{equ}
  G_{y,k}(\tilde x,\tilde y,t) = \cases{G(\tilde x,\tilde y)&if $V(\tilde y) \le k\bigl(V(y) + t^2\bigr)$,\cr 0&otherwise.}
\end{equ}
As above, we define
\begin{equ}
   \Dens^{n,k}_{x,y}(\omega) = \exp\Bigl(\int_0^n G_{y,k}\bigl(\cev\Phi_t(x,y,\omega),\Phi_t(y,\omega),t\bigr)\,d\omega(t) - {1\over 2}\int_0^n \|G_{y,k}(\ldots)\|^2\,dt\Bigr)\;.
\end{equ}
By definition, $\Dens^{n,k}_{x,y}(\omega) = \Dens^{n}_{x,y}(\omega)$ for every $\omega\in B_{y,k}$. Setting $\xi = \delta(\alpha+\beta)$, we thus have the estimate
\begin{equs}
 \int_{\!A_{y,k}} &\bigl( \Dens^{n}_{x,y}(\omega)\bigr)^{-2}\,\PP^n(d\omega)
\le \int_{\!B_{y,k}} \bigl( \Dens^{n,k}_{x,y}(\omega)\bigr)^{-2}\,\PP^n(d\omega) \\
&\le \biggl(\int_{\!B_{y,k}} \exp\Bigl(10 \int_0^n \|G_{y,k}\bigl(\cev\Phi_t(x,y,\omega),\Phi_t(y,\omega),t\bigr)\|^2\,dt\Bigr)\,\PP^n(d\omega)\biggr)^{1/2}\\
&\le \biggl(\int_{\!B_{y,k}} \exp\biggl(10 \int_0^n Ce^{-\gamma t}\Bigl(1+ V(x) + \sup_{s\le t} V\bigl(\Phi_s(y,\omega)\bigr)\Bigr)^\xi\,dt\biggr)\,\PP^n(d\omega)\biggr)^{1/2}\\
&\le \exp\Bigl(C\!\int_0^n e^{-\gamma t}\bigl(1 + k\tilde V(x,y) + kt^2\bigr)^\xi\,dt\Bigr)\;.
\end{equs}
In this expression, we used the Cauchy-Schwarz inequality to go from the first to the second line, and we used assumptions~\ref{ass:ExpSmall} and \ref{ass:Conv2} to go from the second to the third line. Since the integral converges for $n\to\infty$, the bound is uniform in $n$ and Assumption~\ref{ass:Prob} is verified.

The verification of Assumption~\ref{ass:Coupl} is quite similar. Fix some positive constant $K > 0$ and use again the cutoff function
\begin{equ}
  \GN(\tilde x,\tilde y) = \cases{G(\tilde x, \tilde y)&if $\|G(\tilde x,\tilde y)\|^2 \le N$,\cr 0&otherwise.}
\end{equ}
The precise value of $N$ (as a function of $K$) will be fixed later.
We also fix a pair $(x_0,y_0)\in\CH^2$ with $\tilde V(x_0,y_0) \le K$, a value $n > 0$, and initial conditions $(x,y) \in Q^n_K(x_0,y_0)$. By the definition of $Q^n_K(x_0,y_0)$, there exists an element $\tilde \omega \in\SP^n$ such that
\begin{equ}[e:equxy]
(x,y) = \bigl(\cev\Phi_n(x_0,y_0,\tilde\omega),\Phi_n(y_0,\tilde\omega)\bigr)\;,
\end{equ}
and such that
\begin{equ}[e:propw]
\sup_{s\in[0,n]} \tilde V\bigl(\cev\Phi_s(x_0,y_0,\tilde\omega),\Phi_s(y_0,\tilde\omega)\bigr) \le K\;.
\end{equ}
Following the statement of Assumption~\ref{ass:Coupl}, we define the set
\begin{equ}
 B^K_{x,y} = \Bigl\{\omega\in\SP\;\Big|\; \sup_{t\in[0,1]} \tilde V\bigl(\cev\Phi_t(x,y,\omega), \Phi_t(y,\omega)\bigr) \le K\Bigr\}\;,
\end{equ}
which is equal in our setup to the set over which integration is performed in \eref{e:ass4}.
Being now accustomed to these truncation procedures, we define again
\begin{equ}
    \Dens^{(K)}_{x,y}(\omega) = \exp\Bigl(\int_0^1 \GN\bigl(\cev\Phi_t(x,y,\omega),\Phi_t(y,\omega)\bigr)\,d\omega(t) - {1\over 2}\int_0^1 \|\GN(\ldots)\|^2\,dt\Bigr)\;.
\end{equ}
By \eref{e:equxy} and the cocycle property, we can write the integral in
the above expression as
\begin{equ}
 \int_0^{1} \GN\bigl(\cev\Phi_{n+t}(x_0,y_0,\tilde\omega\omega),\Phi_{n+t}(y_0,\tilde\omega\omega)\bigr)\,d\omega(t)\;,
\end{equ}
where $\tilde\omega\omega$ is the realization of the noise which is equal to $\tilde\omega$ for a time $n$ and then to $\omega$ for a time $1$.
Using the {\em a-priori} bound \eref{e:propw} as well as assumptions~\ref{ass:ExpSmall} and \ref{ass:Conv2}, we thus see that there exists a constant $C$ such that the choice
\begin{equ}
 N = C e^{-\alpha \gamma n} \bigl(1 + K\bigr)^{\alpha + \beta}\;,
\end{equ}
ensures that $\Dens^{(K)}_{x,y}(\omega)$ is equal to $\Dens_{x,y}(\omega)$ for $\omega\in B^K_{x,y}$.

We then have the estimate
\begin{equs}
 \int_{B^K_{x,y}}&\bigl(1 - \Dens_{x,y}(\omega)\bigr)^2\,\PP(d\omega)
\le \int_{\SP}\bigl(1 - \Dens^{(K)}_{x,y}(\omega)\bigr)^2\,\PP(d\omega)
\\
&= \int_{\SP}\bigl(\Dens^{(K)}_{x,y}(\omega)\bigr)^2\,\PP(d\omega) - 1 \\
&\le  \biggl(\int_{\SP} \exp\Bigl(6\int_0^1 \|\GN\bigl(\cev\Phi_t(x,y,\omega),\Phi_t(y,\omega)\bigr)\|^2\,dt\Bigr)\,\PP(d\omega)\biggr)^{1/2} - 1\\
&\le  \exp\bigl(C e^{-\alpha \gamma n}(1+K)^{\alpha + \beta}\bigr)-1\;.
\end{equs}
If we take $n \ge \beta \ln(1+K)/\gamma$, the exponent is bounded by $C$ and there exists a constant $C'$ such that
\begin{equ}
 \int_{B^K_{x,y}}\bigl(1 - \Dens_{x,y}(\omega)\bigr)^2\,\PP(d\omega)
\le C' e^{-\alpha \gamma n}(1+K)^{\alpha + \beta}\;,
\end{equ}
thus validating Assumption~\ref{ass:Coupl} with $\gamma_2 = \alpha\gamma$ and $\zeta = \alpha + \beta$.

The proof of \theo{theo:mainEqu} is complete.
\end{proof}

\section{Examples}
\label{sect:examples}

Numerous recent results show that the invariant measure for the 2D Navier-Stokes equation (and also for other dissipative PDEs) is unique if a sufficient number of low-frequency modes are forced by the noise \cite{BC,BKLExp,EMS,MatNS,ELiu,KS00,ArmenKuk,MY}. These results are not covered directly by \theo{theo:mainEqu}, but some more work is needed. The reason is that the sets $A_x^k$ defined in \eref{e:defAxk} are not the natural sets that allow to control the influence of the low-frequency modes onto the high-frequency modes in the 2D Navier-Stokes equation.

On the other hand, our formulation of \theo{theo:mainEqu} makes it quite
easy to verify that the $n$-dimensional Ginzburg-Landau equation (in a bounded domain) shows exponential mixing properties, if sufficiently many low-frequency modes are forced by the noise. We verify this in the following subsection.

\subsection{The Ginzburg-Landau equation}

We consider the SPDE given by
\begin{equ}[e:SGL]
 du = \bigl(\Delta u + u - u^3\bigr)\,dt + Q\,dw(t)\;,\qquad u(0) = u_0\;,
\end{equ}
where the function $u$ belongs to the Hilbert space
\begin{equ}
 \CH = \Ltwo\bigl([-L,L]^n,\R\bigr)\;,
\end{equ}
and $\Delta$ denotes the Laplacean with periodic boundary conditions. The symbol $Q\,d\omega(t)$ stands as a shorthand for
\begin{equ}
 Q\,d\omega(t) \equiv \sum_{i=1}^N q_i e_i\,d\omega_i(t)\;,
\end{equ}
where $\{q_i\}_{i=1}^N$ is a collection of strictly positive numbers, $e_i$ denotes the $i$th eigenvector of the Laplacean, and the $\omega_i$ are $N$ independent Brownian motions (for some finite integer $N$). We also denote by $\lambda_i$ the eigenvalue of $\Delta$ corresponding to $e_i$ and we assume that they are ordered by $\ldots \le \lambda_2\le \lambda_1 \le 0$. We will see that it is fairly easy to construct a binding function $G$ for which the assumptions of the previous section hold with $V(u) = \|u\|$, where $\|\cdot\|$ denotes the norm of $\CH$.

In \cite{ZDP}, it is shown that \eref{e:SGL} possesses a unique mild solution for initial conditions $u_0\in\L^\infty([-L,L]^n)$. It is straightforward to extend this to every initial condition $u_0\in\CH$, by using the regularizing properties of the heat semigroup. Thus, Assumption~\ref{ass:ExistSol} holds and we denote by $\pP^t_u$ the transition probabilities of the solution at time $t$ starting from $u$. We have the following result:

\begin{theorem}\label{theo:SGL}
There exist positive constants $C$ and $\gamma$, and a unique measure $\mu_*\in\Meas_1(\CH)$ such that
\begin{equ}[e:estGL]
 \|\pP^t_u-\mu_*\|_\Lip \le Ce^{-\gamma t}(1 + \|u\|)\;,
\end{equ}
for every $u\in\CH$ and every $t>0$.
\end{theorem}

\begin{proof}
We verify that the assumptions of \theo{theo:mainEqu} hold.
The bounds required for the verification of Assumption~\ref{ass:Lyap2} can be found in \cite{Cerr,ZDP}, for example.

It remains to construct the forcing $G\from\CH^2\to\R^N$ and to verify assumptions~\ref{ass:ExpSmall} and \ref{ass:Conv2}. We consider two copies $u_1$ and $u_2$ of \eref{e:SGL}, with the noise $d\omega$ replaced by $d\omega + G\,dt$ in the second copy. We also denote by $\rho = u_2 - u_1$ the difference process. It satisfies the differential equation
\begin{equ}[e:defrhoGL]
 \dot\rho = \Delta \rho + \rho - \rho\bigl(u_1^2 + u_1 u_2 + u_2^2\bigr)\ + Q\,G(u_1,u_2)\;.
\end{equ}
We can project \eref{e:defrhoGL} onto the direction given by $e_k$. This yields
\begin{equ}
 \dot \rho_k = (\lambda_k+1) \rho_k - \Bigl(\rho\bigl(u_1^2 + u_1 u_2 + u_2^2\bigr)\Bigr)_k + q_k G_k(u_1,u_2)\;,
\end{equ}
for $k=1,\ldots,N$ and
\begin{equ}
  \dot \rho_k = (\lambda_k+1) \rho_k - \Bigl(\rho\bigl(u_1^2 + u_1 u_2 + u_2^2\bigr)\Bigr)_k\;,
\end{equ}
for $k>N$.
We choose $G_k$ for $k=1,\ldots,N$ as
\begin{equ}
 G_k(u_1,u_2) = - {2+\lambda_k\over q_k}\rho_k\;.
\end{equ}
Since $G_k$ can only be defined this way if $q_k\neq 0$, we use at this point the fact that the noise acts directly and independently on {\it every unstable mode}. This requirement can be significantly weakened with the help of \theo{theo:mainEqu}.
We will focus next on more degenerate problems which illustrate the power of our technique.

This choice satisfies Assumption~\ref{ass:ExpSmall}.
With this choice, we can write down the evolution of the norm of $\rho$ as
\begin{equ}
 {d\|\rho\|^2\over dt} = 2\scal{\rho,A\rho} - 2\scal[b]{\rho,\rho\bigl(u_1^2 + u_1 u_2 + u_2^2\bigr)}\;,
\end{equ}
with $A$ the linear operator given by adding up the contribution of $\Delta + 1$ and the contribution of $G$.
By the condition we imposed on $N$, there exists a constant $a>0$ such that $\scal{\rho,A\rho} \le -a\|\rho\|^2$. Furthermore, one has
\begin{equ}
 \scal[b]{\rho,\rho\bigl(u_1^2 + u_1 u_2 + u_2^2\bigr)} \ge 0\;.
\end{equ}
We thus have the differential inequality
\begin{equ}
  {d\|\rho\|^2\over dt} \le -2a \|\rho\|^2 \;,
\end{equ}
which implies that
\begin{equ}
 \|\rho(t)\| \le e^{-at} \|\rho(0)\|\;.
\end{equ}
This implies by \rem{rem:boundV} that Assumption~\ref{ass:Conv2} is also satisfied. The proof of \theo{theo:SGL} is complete.
\end{proof}

\subsection{A reaction-diffusion system}

Consider the following reaction-diffusion system:
\begin{equa}[e:toy2]
 du &= \bigl(\Delta u + 2u + v - u^3\bigr)\,dt + dw(t)\;,\\
 dv &= \bigl(\Delta v + 2v + u - v^3\bigr)\,dt\;,
\end{equa}
where the pair $(u,v)$ belongs to the Hilbert space
\begin{equ}
 \CH = \CH_u \oplus \CH_v = \Ltwo\bigl([-L,L],\R\bigr)\oplus \Ltwo\bigl([-L,L],\R\bigr)\;.
\end{equ}
The symbol $\Delta$ again denotes the Laplacean with periodic boundary conditions,
and $d\omega$ is the cylindrical Wiener process on $\CH_u$ (meaning that it is space-time white noise).

Notice that, because of the presence of $v$, this system does not satisfy the assumptions stated in the papers mentioned at the beginning of this section. In other words, even though the forcing is infinite-dimensional, not all the determining modes for \eref{e:toy2} are forced.

We take as our Lyapunov function
\begin{equ}
 V(u,v) = \|u\|_\infty + \|v\|_\infty\;,
\end{equ}
with $\|\cdot\|_\infty$ the $\L^\infty$ norm. As in the previous subsection, one can show that with this choice of $V$, our problem satisfies assumptions~\ref{ass:ExistSol} and \ref{ass:Lyap2}.
We will now construct a binding function $G$ which satisfies assumptions~\ref{ass:Conv2} and \ref{ass:ExpSmall}. We consider, as in \eref{e:defcoupEqu}, two copies $(u_1,v_1)$ and $(u_2,v_2)$ of the system \eref{e:toy2}, but the noise is modified by $G$ on the second copy. We also define $\rho_u = u_2 - u_1$ and $\rho_v = v_2 - v_1$. We then have
\begin{equa}[e:defrhouv]
 \dot \rho_u &= \Delta \rho_u + 2\rho_u + \rho_v - \rho_u (u_1^2 + u_1 u_2 + u_2^2) + G(u_1,u_2,v_1,v_2)\;,\qquad\\
 \dot \rho_v &= \Delta \rho_v + 2\rho_v + \rho_u - \rho_v (v_1^2 + v_1 v_2 + v_2^2)\;.
\end{equa}
Our construction of $G$ is inspired from the construction we presented in \sect{sect:toy}. We introduce the variable $\zeta = \rho_u + 3\rho_v$. Substituting this in \eref{e:defrhouv}, it defines the function $G$ if we impose that the equation for $\zeta$ becomes
\begin{equ}[e:zeta]
 \dot \zeta = \Delta\zeta - \zeta\;,
\end{equ}
so that $\|\zeta(t)\|^2 \le \|\zeta(0)\|^2 e^{-t}$.
Notice that the function $G$ achieving this identity satisfies a bound of the type
\begin{equ}
 \|G\| \le C\bigl(\|\rho_u\| + \|\rho_v\|\bigr)\bigl(1 + \|u_1\|_\infty +  \|u_2\|_\infty + \|v_1\|_\infty + \|v_2\|_\infty\bigr)^2\;,
\end{equ}
thus satisfying Assumption~\ref{ass:ExpSmall}. It remains to show that Assumption~\ref{ass:Conv2} is satisfied.
The equation for $\rho_v$ reads
\begin{equ}
  \dot \rho_v = \Delta \rho_v - \rho_v + \zeta - \rho_v (v_1^2 + v_1 v_2 + v_2^2)\;.
\end{equ}
Therefore, the norm of $\rho_v$ satisfies
\begin{equ}
 \|\rho_v(t)\|^2 \le \|\rho_v(0)\|^2e^{-t} + {1+\|\zeta(0)\|^2\over 2}e^{-t}\;.
\end{equ}
This in turn implies, through the definition of $\zeta$ and the bound on $\|\zeta(t)\|$, that a similar bound holds for $\|\rho_u(t)\|$.
This shows that the bound \eref{e:boundrho} is satisfied.
Similar estimates hold with the $\L^\infty$ norm replacing the $\L^2$ norm, and so Assumption~\ref{ass:Conv2} is satisfied by \rem{rem:boundV}.

In fact, a straightforward computation, which can be found in \cite{Cerr,H01,MSUnif} for example, shows that in this example, one can get a uniform estimate on the Lyapunov function $V$. More precisely, there exists a constant $C$ such that for all initial conditions $x\in\CH$,
\begin{equ}[e:estV]
 \int_{\CH} V(y)\,\pP^_x(dy) \le C\;.
\end{equ}
Denoting by $\CP_t^*$ the semigroup acting on measures generated by the solutions of \eref{e:toy2}, we thus have:
\begin{theorem}
There exists a unique probability measure $\mu_* \in \Meas_1(\CH)$ such that $\CP_t^*\mu_* = \mu_*$ for every $t \ge 0$. Furthermore, there exist constants $C$ and $\gamma$ such that
\begin{equ}[e:uniform]
 \|\CP_t^*\nu - \mu_*\|_\Lip \le Ce^{-\gamma t}\;,
\end{equ}
for every $\nu  \in \Meas_1(\CH)$.
\end{theorem}
\begin{proof}
Combining \eref{e:estV} with the results of \theo{theo:main} and a computation similar to what was done in the proof of \cor{cor:main}, we get \eref{e:uniform} for integer times. The generalization to arbitrary times is straightforward, using the fact that the growth rate of the difference process $(\rho_u,\rho_v)$ (with $G \equiv 0$) can easily be controlled.
\end{proof}

\begin{remark}
In fact, the dependence on $u$ in the right-hand side of \eref{e:estGL} can be removed similarly by checking that an estimate of the type \eref{e:estV} is verified for the solutions of the stochastic Ginzburg-Landau equation \eref{e:SGL}.
\end{remark}

\subsection{A chain with nearest-neighbour interactions}

In the previous example, the noise acted on infinitely many degrees of freedom in a non-degenerate way. As a consequence, one step was sufficient to transmit the noise to the entire system. We will now look at a much more degenerate system, where the noise acts on only {\em one} degree of freedom, although an {\em arbitrary} number of modes are linearly unstable.

Our model is given by
\begin{equa}[e:chain]
 dx_0 &= (a^2 x_0 + x_1 - x_0^3)\,dt + d\omega\;,\\
\dot x_k &= (a^2 - k^2)x_k + x_{k-1} + x_{k+1} - x_k^3\;,\qquad k=1,2,\ldots\;,
\end{equa}
where $a\in\R$ is an arbitrary constant. One should think of the deterministic part of \eref{e:chain} as a very simple model for a dissipative PDE of the Ginzburg-Landau type. We will consider \eref{e:chain} in the (real) Hilbert space $\CH=\ell^2$ endowed with its canonical orthonormal basis $\{e_k\}_{k=0}^\infty$. It is easy to verify that \eref{e:chain} possesses a unique solution. We denote again by $\CP_t^*$ the semigroup acting on measures $\nu\in\Meas(\ell^2)$ generated by \eref{e:chain}. We will show

\begin{theorem}\label{theo:chain}
For the problem \eref{e:chain}, there exists a unique probability measure $\mu_* \in \Meas_1(\ell^2)$ such that $\CP_t^*\mu_* = \mu_*$ for every $t \ge 0$. Furthermore, there exist constants $C$ and $\gamma$ such that
\begin{equ}
 \|\CP_t^*\nu - \mu_*\|_\Lip \le Ce^{-\gamma t}\;,
\end{equ}
for every $\nu  \in \Meas_1(\CH)$.
\end{theorem}

\begin{proof}
We will take as our Lyapunov function $V(x) = \|x\|^p$ for some power of $p$ to be fixed later. It is a straightforward task to verify that the dynamics generated by \eref{e:chain} does indeed satisfy assumptions~\ref{ass:ExistSol} and \ref{ass:Lyap2} for this choice of $V$. 

We next show that a bound of the type \eref{e:estV} holds for the solutions of \eref{e:chain}, thus yielding the uniformity in the convergence towards the invariant measure $\mu_*$. 
Let us define the process $y(t) \in\ell^2$ by $y(t) = x(t) - \omega(t) e_0$. This process then satisfies the following system of differential equations:
\begin{equa}[e:chainy]
 \dot y_0 &= a^2 (y_0 + \omega) + y_1 - (y_0+ \omega)^3\;,\\
 \dot y_1 &= (a^2-1) y_1 + y_0 + y_2 - y_1^3 + \omega\;,\\
\dot y_k &= (a^2 - k^2)y_k + y_{k-1} + y_{k+1} - y_k^3\;,\qquad k=2,3,\ldots\;.
\end{equa}
We denote by $\|y\|_\infty$ the norm given by $\sup_k |y_k|$. It follows from \cite{Lun} that \eref{e:chainy} possesses a strong solution for positive times. Furthermore, from \eref{e:chainy} and the definition of the $\|\cdot\|_\infty$-norm, we see that there are constants $c_1,c_2>0$ such that
\begin{equ}[e:estinfty]
{D_-\|y\|_\infty \over Dt} \le - c_1 \|y\|_\infty^3 + c_2 \bigl(1 + |\omega(t)|^3\bigr)\;,
\end{equ}
where $D_-/Dt$ denotes the left-handed lower Dini derivative. A straightforward computation shows that \eref{e:estinfty} implies that there exists a constant $C$ such such that
\begin{equ}
 \|y(1/2)\|_\infty \le C\sup_{t\in[0,1/2]} \bigl(1 + |\omega(t)|\bigr)\;,
\end{equ}
independently of the initial condition. In order to conclude the proof of the estimate \eref{e:estV}, it suffices to show that there exists a constant $C$ such that
\begin{equ}
 \Exp \|y(1/2)\| \le C\bigl(1+\|y(0)\|_\infty\bigr)\;.
\end{equ}
This follows easily from the dissipativity of the nonlinearity in $\CH$ and the fact that the semigroup generated by the linear part of \eref{e:chainy} is bounded from $\ell^\infty$ into $\ell^2$.

It remains to verify that the assumptions \ref{ass:ExistSol}--\ref{ass:Conv2} are indeed satisfied for some binding function $G$. This, together with the uniform bound obtained above, shows that the conclusions of \theo{theo:chain} hold.
As for the toy model presented in \sect{sect:toy}, we consider a process $y \in \ell^2$ governed by the same equation as \eref{e:chain}, but with $d\omega$ replaced by $d\omega + G(x,y)\,dt$. We then introduce the difference process $\rho = y-x$, which is given by the solution of
\begin{equa}[e:rhochain]
 \dot \rho_0 &= a^2 \rho_0 + \rho_1 - \rho_0(x_0^2 + x_0 y_0 + y_0^2) + G(x,y)\;,\\
 \dot \rho_k &= (a^2 - k^2) \rho_k + \rho_{k+1} + \rho_{k-1} - \rho_k(x_k^2 + x_k y_k + y_k^2)\;.
\end{equa}
The aim of the game is to find a function $G$ for which $\rho(t)\to 0$ as $t\to\infty$. We can split \eref{e:rhochain} into ``low modes'' and ``high modes'' by introducing
\begin{equ}
 k_* = \inf\{k>0\;|\; k^2 - a^2 \ge 3\}\;.
\end{equ}
At the level of the Hilbert space $\ell^2$, we set $\ell^2 = \CH_L \oplus \CH_H$, where $\CH_L\approx \R^{k_*}$ is generated by $e_0,\ldots,e_{k_*-1}$ and $\CH_H$ is its orthogonal complement. We denote by $\rho_L$ and $\rho_H$ the components of $\rho$ and by $A_H$ the restriction (as a symmetric quadratic form) of the linear part of \eref{e:chain} to $\CH_H$. It is by construction easy to see that
\begin{equ}
 \scal{\rho_H,A_H\rho_H} \ge \|\rho_H\|^2\;.
\end{equ}
As a consequence, we have for $\|\rho_H\|^2$ the following estimate:
\begin{equ}[e:estrhoH]
 \|\rho_H(t)\|^2 \le e^{-t}\|\rho_H(0)\|^2 + {1\over 4}\int_0^t e^{t-s}|\zeta_1(s)|^2\,ds\;,
\end{equ}
where we defined $\zeta_1 = \rho_{k_*-1}$. (The reason for renaming $\rho_{k_*-1}$ this way will become clear immediately.)
It remains to construct $G$ in such a way to get good estimates on $\|\rho_L(t)\|^2$. In order to achieve this, we use again the same method as for the first toy model. The variable $\zeta_1$ obeys the equation
\begin{equ}
\dot \zeta_1 = c_1 \rho_{k_*-1} + \rho_{k_*} + \rho_{k_*-2} - \rho_{k_*-1}(x_{k_*-1}^2 + x_{k_*-1} y_{k_*-1} + y_{k_*-1}^2)\;,
\end{equ}
with some constant $c_1\in\R$. We thus introduce a new variable $\zeta_2$ defined by
\begin{equ}
 \zeta_2 = (c_1+1) \rho_{k_*-1} + \rho_{k_*} + \rho_{k_*-2} - \rho_{k_*-1}(x_{k_*-1}^2 + x_{k_*-1} y_{k_*-1} + y_{k_*-1}^2)\;.
\end{equ}
It is important to notice two facts about this definition. The first is that it yields for $|\zeta_1|^2$ the estimate
\begin{equ}[e:est1]
  |\zeta_1(t)|^2 \le e^{-t}|\zeta_1(0)|^2 + {1\over 4}\int_0^t e^{t-s}|\zeta_2(s)|^2\,ds\;.
\end{equ}
The second is that $\zeta_2$ can be written in the form
\begin{equ}
 \zeta_2 = \rho_{k_*-2} + \CQ_2(\rho,x,y)\;,
\end{equ}
where $\CQ_2$ is a polynomial depending only on components $\rho_i$, $x_i$ and $y_i$ with $i \ge k_*-1$, and such that each of its terms contains at least one factor $\rho_i$.

Now look at the equation for $\dot \zeta_2$. It is clear from the structure of $\zeta_2$ and from the structure of the equations \eref{e:chain} and \eref{e:rhochain} that it can be written as
\begin{equ}
 \dot \zeta_2 = -\zeta_2 + \zeta_3\;,
\end{equ}
where
\begin{equ}
 \zeta_3 = \rho_{k_*-3} + \CQ_3(\rho,x,y)\;.
\end{equ}
This time, the polynomial $\CQ_3$ depends only on components with an index $i \ge k_*-2$. This procedure can be iterated, yielding a whole family of variables
\begin{equ}[e:defzetal]
 \zeta_l = \rho_{k_*-l} + \CQ_l(\rho,x,y)\;,
\end{equ}
where the $\CQ_l$ are polynomials depending only on indices $i\ge k_*-l+1$, and containing at least one factor $\rho_i$ in each term. Furthermore, one gets for every $\zeta_l$ the estimate
\begin{equ}[e:estzetal]
   |\zeta_l(t)|^2 \le e^{-t}|\zeta_l(0)|^2 + {1\over 4}\int_0^t e^{t-s}|\zeta_{l+1}(s)|^2\,ds\;.
\end{equ}
Notice that \eref{e:estzetal} is valid for $l < k_*$. For $l = k_*$, we have
\begin{equ}[e:lastequ]
\dot \zeta_{k_*} = \CQ_{k_*+1}(\rho,x,y) + G(x,y)\;.
\end{equ}
It thus suffices to choose $G$ in such a way that \eref{e:lastequ} becomes
\begin{equ}[e:estzetalast]
 \dot \zeta_{k_*} = -\zeta_{k_*}\;.
\end{equ}
Denoting by $\zeta$ the vector $\zeta_1,\ldots,\zeta_{k_*}$, we get from
\eref{e:estzetal} and \eref{e:estzetalast} the estimate
\begin{equ}[e:estzeta]
 \|\zeta(t)\|^2 \le C e^{-\gamma t}\|\zeta(0)\|^2\;,
\end{equ}
for any $\gamma \in (0,1)$. Plugging this into \eref{e:estrhoH} yields for $\|\rho_H\|$ the estimate
\begin{equs}
  \|\rho_H(t)\|^2 &\le C e^{-\gamma t}\bigl(\|\rho_H(0)\|^2 + \|\zeta(0)\|^2\bigr)\\
&\le C e^{-\gamma t}\|\rho(0)\|^2\bigl(1 + \|x(0)\| + \|y(0)\|\bigr)^p\;,
\end{equs}
for some constants $C$, $\gamma$ and $p$. It remains to get an estimate on $\|\rho_L\|$. From \eref{e:estzeta} and the definition of $\zeta_1$, we get immediately
\begin{equ}
 |\rho_{k_*-1}(t)|^2\le C e^{-\gamma t}\|\rho(0)\|^2\bigl(1 + \|x(0)\| + \|y(0)\|\bigr)^p\;.
\end{equ}
From the definition of $\zeta_2$, we get
\begin{equ}
 |\rho_{k_*-2}(t)|^2\le C\Bigl(|\zeta_2(t)|^2 + \bigl|\CQ_2\bigl(\rho(t),x(t),y(t)\bigr)\bigr|^2\Bigr)\;.
\end{equ}
But we know that $\CQ_2$ only depends on components of $\rho$, $x$, and $y$ with an index $i \ge k_*-1$. These are precisely the components of $\rho$ on which we already have an estimate. We thus get
\begin{equ}
  |\rho_{k_*-2}(t)|^2\le C e^{-\gamma t}\|\rho(0)\|^2\bigl(1 + \|x(0)\| + \|y(0)\| + \|x(t)\|\bigr)^p\;,
\end{equ}
for some other power $p$. Here we used the fact that $y(t)=x(t)+\rho(t)$ to get rid of $\|y(t)\|$ in the estimate. The same reasoning can be applied to $\rho_{k_*-3}$, and so forth down to $\rho_0$. We finally get
\begin{equ}[e:estlow]
   \|\rho_L(t)\|^2\le C e^{-\gamma t}\|\rho(0)\|^2\bigl(1 + \|x(0)\| + \|y(0)\| + \|x(t)\|\bigr)^p\;,
\end{equ}
for some (large) power of $p$. We thus verified \eref{e:boundrhodir}. The bound \eref{e:boundrhoinv} is obtained in the same way, by noticing that we can as well get the estimate
\begin{equ}
   \|\rho_L(t)\|^2\le C e^{-\gamma t}\|\rho(0)\|^2\bigl(1 + \|x(0)\| + \|y(0)\| + \|y(t)\|\bigr)^p\;,
\end{equ}
instead of \eref{e:estlow}. The proof of \theo{theo:chain} is complete.
\end{proof}

\begin{remark}
The whole construction is strongly reminiscent of what was done in \cite{EPR2} to control a finite Hamiltonian chain of non-linear oscillators with nearest-neighbour coupling driven by thermal noise at its boundaries.
\end{remark}

\begin{remark}
The linearity of the nearest-neighbour coupling is not essential for our argument. We could as well have replaced \eref{e:chain} by
\begin{equa}
 dx_0 &= \bigl(a^2 x_0 + V_2'(x_1-x_0) - V_1'(x_0)\bigr)\,dt + d\omega\;,\\
\dot x_k &= (a^2 - k^2)x_k + V_2'(x_{k-1}- x_k) + V_2'(x_{k+1}-x_k) - V_1'(x_k)\;,
\end{equa}
with $V_1$ and $V_2$ two polynomial-like functions, \ie smooth functions such that
\begin{equ}
 {d^n V_i(x)\over dx^n} \approx x^{\alpha_i-n} \qquad\text{for}\qquad
|x| \to \infty\;,
\end{equ}
for some $\alpha_i \ge 2$.
Imposing the condition $V_2''(x) \ge c$ for some $c>0$ yields an effective coupling between neighbours at every point of the phase space. This is sufficient to apply our construction.
\end{remark}

\bibliographystyle{myalph}
\markboth{\sc \refname}{\sc \refname}
\bibliography{refs}

\end{document}